\documentclass{amsart}
\usepackage[margin=1in]{geometry}
\usepackage{amsmath, amsfonts, amssymb}
\usepackage{subfig}
\usepackage{graphicx}
\usepackage{color}
\allowdisplaybreaks

\graphicspath{{images/}}

\title{A surface hopping Gaussian beam method for high-dimensional transport systems}

\author{Zhenning Cai}
\address[Zhenning Cai]{Department of Mathematics, National University of Singapore,
  Level 4, Block S17, 10 Lower Kent Ridge Road, Singapore 119076}
\email{matcz@nus.edu.sg}

\author{Jianfeng Lu}
\address[Jianfeng Lu]{Department of Mathematics, Department of Physics, Department of Chemistry,
  Duke University, Box 90320, Durham NC 27708, USA}
\email{jianfeng@math.duke.edu}

\thanks{This work is partially supported by the National Science
  Foundation under Grant Nos. DMS-1454939 and RNMS11-07444
  (KI-Net). The authors would like to thank Giovanni Ciccotti for
  illuminating discussions.}

\DeclareMathOperator{\re}{Re}
\DeclareMathOperator{\im}{Im}

\newtheorem{proposition}{Proposition}
{\theoremstyle{remark} \newtheorem*{remark*}{Remark}}

\begin{document}
\maketitle

\begin{abstract}
  We propose a surface hopping Gaussian beam method to numerically
  solve a class of high frequency linear transport systems in high
  spatial dimensions, based on asymptotic analysis. The stochastic
  surface hopping is combined with Gaussian beam method to deal with
  the multiple characteristic directions of the transport system in
  high dimensions. The Monte Carlo nature of the proposed algorithm
  makes it easy for parallel implementations. We validate the
  performance of the algorithms for applications on the
  quantum-classical Liouville equations.
\end{abstract}

\section{Introduction}
In this paper we focus on the system of $u_i^{\epsilon}$,
$i=1,\cdots,n$ with the following form:
\begin{align} \label{eq:equations}
& \frac{\partial u_i^{\epsilon}(t,x)}{\partial t} +
  \alpha_i(t,x)^{\top} \nabla_x u_i^{\epsilon}(t,x) =
\frac{\mathrm{i}}{\epsilon} \beta_i(t,x) u_i^{\epsilon}(t,x) +
  \sum_{j=1}^n \gamma_{ij}(t,x) u_j^{\epsilon}(t,x),  \\
& \label{eq:init}
u_i^{\epsilon}(0,x) = v_i^{\epsilon}(x), \qquad i=1,\cdots,n.
\end{align}
where $\epsilon$ is a small positive parameter and $v_i^{\epsilon}$
are initial conditions. For each $i$, the unknown
$u_i^{\epsilon}(\cdot,\cdot)$ maps from
$\mathbb{R^+} \times \mathbb{R}^m$ to $\mathbb{C}$, where $m$ is large
so that this is a system in high spatial dimension. For each $i, j$,
The coefficients $\alpha_i(\cdot,\cdot)$ maps from
$\mathbb{R^+} \times \mathbb{R}^m$ to $\mathbb{R}^m$,
$\beta_i(\cdot, \cdot)$ maps from $\mathbb{R^+} \times \mathbb{R}^m$
to $\mathbb{R}$, and $\gamma_{ij}(\cdot,\cdot)$ maps from
$\mathbb{R^+} \times \mathbb{R}^m$ to $\mathbb{C}$.
If $n=1$, the equation \eqref{eq:equations} becomes a linear
transport equation with a stiff source term. When $n > 1$, the system
consists of $n$ such ingredients and they interact with each other
through the second linear term on the right hand side. For simplicity,
here we assume that all the coefficients $\alpha_i(\cdot,\cdot)$,
$\beta_i(\cdot,\cdot)$ and $\gamma_{ij}(\cdot,\cdot)$ are smooth
functions.

Systems with the form \eqref{eq:equations} arise in several
application fields.  In this work, we focus in particular on the
quantum-classical Liouville equation (QCLE) \cite{Kapral1999,
  Schutte1999}, which is a model for mixed quantum-classical dynamics
(say for quantum systems consist of nuclei and electrons). The
equations are derived through an adiabatic representation of the
Hamiltonian and a partial Wigner transform of the von Neumann
equation. As a result, the leading orders of the equation have the
form \eqref{eq:equations}, and the terms with $\alpha_i$, $\beta_i$
and $\gamma_{ij}$ in \eqref{eq:equations} respectively represent the
transport of particles, the phase oscillation, and the exchange
between phases \cite{Horenko2002}.  Numerical study of single and dual
crossing examples in \cite{Horenko2002} shows that QCLE can accurately
predict the population of each adiabatic states when comparing with
the fully quantum results. We will come back to more details of QCLE
in our numerical experiments.  Another example of \eqref{eq:equations}
worth mentioning is the multi-component transport equations
\cite{Elste2008,Dou2015} with linear interactions with the form
\begin{displaymath}
\frac{\partial f_i(t,r,p)}{\partial t} =
  -\frac{\partial H_i(r,p)}{\partial p}
    \frac{\partial f_i(t,r,p)}{\partial r} +
  \frac{\partial H_i(r,p)}{\partial r}
    \frac{\partial f_i(t,r,p)}{\partial p} +
  \sum_{j=1}^n \left[
    \gamma_{j\rightarrow i} f_j(t,r,p) -
    \gamma_{i\rightarrow j} f_i(t,r,p)
  \right], 
\end{displaymath}
for $i = 1,\cdots,n$, where $f_i(t,r,p)$ and $H_i(r,p)$ are
respectively the distribution function and the Hamiltonian for the
$i$th component, and $\gamma_{i \rightarrow j}$ is the rate of change
from the $i$th component to the $j$th component.

\smallskip 

While vast research has been done in numerical schemes for
advection-reaction systems in low spatial dimensions (see e.g.,
\cite{Hundsdorfer2003} for a review), the high dimensionality poses great
challenges to conventional numerical schemes. In such cases, the
grid-based methods become prohibitively expensive due to the curse of
dimensionality, and we are therefore motivated to develop
particle-type methods for the system.

Besides the high dimensionality, the fast oscillation of the solution 
due to the $\beta$ term on the right hand side and perhaps the initial
data poses additional challenge to the problem. Semiclassical methods
are typically used to deal with such high frequency wave problems,
such as the Gaussian beam method \cite{Heller1981,Ralston1982,
Hagedorn1998} and also more recently the method based on nonlinear
geometric optics \cite{Crouseilles2016}. 

Let us focus on the Gaussian beam method. The basic idea is to first
approximate the initial data of the equation with the superposition of
many Gaussians with small variance, and then propogate each Gaussian
as asymptotic solution to the equation, the solution  is then given by
the superposition of the Gaussians at a given time.  The oscillation
of the solution is naturally built  into the Gaussian ansatz.  For
scalar equations ($n=1$), the Gaussian beam method has been
extensively studied, especially for the classical Liouville equation
\cite{Ma1993, Horenko2002a,Yin2013}. 
Such methods (including the nonlinear geometric optics method) however
cannot be applied to systems as \eqref{eq:equations} due to the
presence of the zeroth-order ``$\gamma$-term'', since this term may cause the
number of the Gaussian beams to be multiplied by $n$ at each time step, and
eventually leads to an exponential increment of the Gaussian beams, which is
obviously infeasible.

To avoid large amount of Gaussian beams, we consider using a Monte
Carlo method to make it feasible. Thus, at each time step, we have
only one Gaussian beam, but the Gaussian beam may contribute to any
component of the solution randomly according to the $\gamma$-term.
In particular, for the QCLE as mentioned earlier, the $\gamma$-term
appears in the equation corresponding to the ``surface hopping'' of
the trajectories \cite{Kapral1999, Horenko2002, Kapral2016}, which
corresponds to the random choice of the component. Such a perspective
is widely used in understanding the quantum-classical dynamics, and it
is often used to understand the fewest switches surface hopping
algorithm for Schr\"odinger equations \cite{Tully1990,Barbatti2011}.
The basic idea of the surface hopping method is to evolve trajectories
on several ``surfaces'', and let the trajectories
stochastically hop from one surface to another from time to time. This
algorithm is mostly applied to the Schr{\"o}dinger equation, and a
recent paper by one of the authors \cite{Lu2016} provides a
mathematical understanding of the fewest switch surface hoppping
algorithm \cite{Tully1990}. The surface hopping method for the
quantum-classical Liouville equation has been considered in
\cite{Horenko2002, MacKernan2002,Kernan2008, Kelly2013}, while those
methods are heuristically proposed without rigorous error control.


The main contribution of our work is a novel
semiclassical approach that generalizes Gaussian beam method for
systems like \eqref{eq:equations}. We call the proposed method the
\emph{surface hopping Gaussian beam method}, as it combines the ideas of
Gaussian beam method and surface hopping to solve
\eqref{eq:equations}. The algorithm allows each trajectory to be
evolved independently, and thus its parallelization becomes
trivial. In order to motivate the algorithm clearly, we first derive
the Gaussian beam part and the surface hopping part of the algorithm
independently, and then combine them together. The proposed method is
verified both theoretically and numerically. This work also depicts
the surface hopping in a clearer and more general sense, especially
from the mathematical point of view.

In the remaining part of this paper, the proposed surface hopping
Gaussian beam method is constructed in Section \ref{sec:method}, where
the Gaussian beam method and the surface hopping method are developed
respectively in Section \ref{sec:gb} and Section \ref{sec:sh}, and
Section \ref{sec:shgb}--\ref{sec:init_data} describe and validate the
whole algorithm with initial data processing. Numerical examples are
given in Section \ref{sec:num_ex}. The whole paper is briefly
summarized in Section \ref{sec:summary}.

\section{Construction of the surface hopping Gaussian beam method}
\label{sec:method}

For simplicity, let us first assume that
the initial data $v_i^{\epsilon}$ in \eqref{eq:init} has the following
form:
\begin{equation}
v_i^{\epsilon}(x) =
  G^{\epsilon}(M_0, N_0, X_0, P_0, S_0, A_0; x) \delta_{ii_0}, \qquad
i=1,\cdots,n,
\end{equation}
where $M_0$ and $N_0$ are $m \times m$ real symmetric matrices, 
$X_0, P_0 \in \mathbb{R}^m$, $S_0 \in \mathbb{R}$,
$A_0 \in \mathbb{C}$, and $i_0 \in \{1,\cdots,n\}$. The function
$G^{\epsilon}$ is given by 
\begin{equation}
G^{\epsilon}(M, N, X, P, S, A; x) = A \exp \left(
  - \frac{(x-X)^{\top} (M + \mathrm{i}N) (x-X)}{2\epsilon}
  + \frac{\mathrm{i} P^{\top} (x-X)}{\epsilon}
  + \mathrm{i} \frac{S}{\epsilon}
\right).
\end{equation}
Apparently $G^{\epsilon}$ defines a ``Gaussian beam'' with width
$\sqrt{\epsilon}$, amplitude $A$, and center $X$. The matrix $M$ is required
to be positive definite so that the function decays to zero at
infinity. General initial conditions will be represented as
superposition of such Gaussian beams as discussed in
Section~\ref{sec:init_data}.

To motivate the surface hopping Gaussian beam method, we will consider
first the following two special cases of the system
\eqref{eq:equations}:
\begin{align}
\label{eq:gb}
& \frac{\partial u_i^{\epsilon}(t,x)}{\partial t} +
  \alpha_i(t,x)^{\top} \nabla_x u_i^{\epsilon}(t,x) =
\frac{\mathrm{i}}{\epsilon} \beta_i(t,x) u_i^{\epsilon}(t,x) +
  \gamma_{ii}(t,x) u_i^{\epsilon}(t,x), \quad i=1,\cdots,n, \\
\intertext{and,}
\label{eq:sh}
& \frac{\partial u_i^{\epsilon}(t,x)}{\partial t} =
  \sum_{\substack{j=1\\ j\neq i}}^n
    \gamma_{ij}(t,x) u_j^{\epsilon}(t,x),
\quad i=1,\cdots,n.
\end{align}
The system \eqref{eq:gb} is obtained by setting all the off-diagonal
entries of the matrix $(\gamma_{ij})$ to be zero in
\eqref{eq:equations}, and the system \eqref{eq:sh} can be obtained by
setting all the coefficients other than the off-diagonal entries of
$(\gamma_{ij})$ to be zero. Note that \eqref{eq:gb} can be viewed as ``diagonal terms'' of \eqref{eq:equations}, while \eqref{eq:sh} captures the ``off-diagonal terms''.

In the remainder of this section, we
will first apply the Gaussian beam method to \eqref{eq:gb}, and then
apply the surface hopping method to \eqref{eq:sh}. The fusion of these
two parts follows thereafter, which eventually becomes a method
solving the system \eqref{eq:equations} in the general case.

\subsection{The Gaussian beam method for the scalar equation}
\label{sec:gb}
The system \eqref{eq:gb} actually consists of $n$ independent scalar
equations for each $u_i^{\epsilon}$, $i=1,\cdots,n$. Since the initial
condition $v_i^{\epsilon}$ is nonzero only if $i = i_0$, we have
$u_i^{\epsilon} \equiv 0$ if $i \neq i_0$, and the system reduces to a
single scalar equation with $i = i_0$ in \eqref{eq:gb}. Here we recall the standard Gaussian beam method applying to the particular case \eqref{eq:gb}. 

In the Gaussian beam method, we take the ansatz 
\begin{equation} \label{eq:ansatz}
u_{i_0}^{\epsilon}(t,x) \approx U_{i_0}^{\epsilon}(t,x) =
  G^{\epsilon}(M(t), N(t), X(t), P(t), S(t), A(t); x).
\end{equation}
By direct calculation, we get the derivatives appearing in
\eqref{eq:gb} as
\begin{align}
\label{eq:time_derivative}
\begin{split}
\frac{\partial U_{i_0}^{\epsilon}(t,x)}{\partial t} &=
  \frac{1}{\epsilon}\bigg[
    (x-X(t))^{\top} [M(t) + \mathrm{i}N(t)]
      \frac{\mathrm{d} X(t)}{\mathrm{d}t}
    - \frac{1}{2} (x-X(t))^{\top} \left(
      \frac{\mathrm{d} M(t)}{\mathrm{d}t} +
      \mathrm{i} \frac{\mathrm{d} N(t)}{\mathrm{d}t}
    \right) (x-X(t)) \\
    & \qquad +\mathrm{i} (x-X(t))^{\top} \frac{\mathrm{d}P(t)}{\mathrm{d}t}
    - \mathrm{i} P(t)^{\top} \frac{\mathrm{d}X(t)}{\mathrm{d}t}
    + \mathrm{i} \frac{\mathrm{d}S(t)}{\mathrm{d}t}
    + \frac{\epsilon}{A(t)} \frac{\mathrm{d}A(t)}{\mathrm{d}t}
  \bigg] U_{i_0}^{\epsilon}(t,x),
\end{split} \\
\label{eq:space_derivative}
\nabla_x U_{i_0}^{\epsilon}(t,x) &= \frac{1}{\epsilon}
\left[ \mathrm{i} P(t) - M(t) (x-X(t)) - \mathrm{i} N(t) (x-X(t)) \right]
  U_{i_0}^{\epsilon}(t,x).
\end{align}
Since
\begin{equation} \label{eq:order}
|x-X|^k G^{\epsilon}(M,N,X,P,S,A;x)
  \sim O(\epsilon^{k/2}), \qquad \forall k \in \mathbb{N},
\end{equation}
for a smooth function $f(t,x)$ independent of $\epsilon$, we have the
following approximation:
\begin{equation} \label{eq:fG}
\begin{split}
f(t,x) G^{\epsilon}(t;x)
&= \bigg[
  f(t, X(t)) + \Big( \nabla_x f(t, X(t)) \Big)^{\top} (x - X(t)) \\
  & \qquad + \frac{1}{2} (x-X(t))^{\top} \nabla_x^2 f(t, X(t)) (x-X(t))
\bigg] G^{\epsilon}(t;x) + O(\epsilon^{3/2}).
\end{split}
\end{equation}
Since we have already assumed the smoothness of $\alpha_{i_0}(t,x)$,
$\beta_{i_0}(t,x)$ and $\gamma_{i_0 j}(t,x)$, the function $f(t,x)$ in
\eqref{eq:fG} can be replaced by these coefficients and the equality
still holds. Thus, by substituting \eqref{eq:time_derivative} and
\eqref{eq:space_derivative} into \eqref{eq:gb}, and applying
the above result, we find that $U_{i_0}^{\epsilon}$ satisfies
\eqref{eq:gb} up to a residual with relative $L^{\infty}$ error
$O(\epsilon^{1/2})$, if the functions $N(t)$, $X(t)$, $P(t)$, $S(t)$
and $A(t)$ fulfill the following equations obtained by matching the
terms with the same orders of $\epsilon$:
\begin{align*}
O(\epsilon^{-1}): & \frac{\mathrm{d}S(t)}{\mathrm{d}t} +
  P(t)^{\top} \left(
    \alpha_{i_0}(t,X(t)) - \frac{\mathrm{d}X(t)}{\mathrm{d}t}
  \right) = \beta_{i_0}(t,X(t)), \\
O(\epsilon^{-\frac{1}{2}}): & 
  M(t) \left(
    \frac{\mathrm{d}X(t)}{\mathrm{d}t} - \alpha_{i_0}(t,X(t))
  \right) = 0, \\
& \frac{\mathrm{d}P(t)}{\mathrm{d}t}
  + [\nabla_x \alpha_{i_0}(t,X(t))]^{\top} P(t)
  + N(t) \left(
    \frac{\mathrm{d}X(t)}{\mathrm{d}t} - \alpha_{i_0}(t,X(t))
  \right) = \nabla_x \beta_{i_0}(t,X(t)), \\
O(1): & \re \left( \frac{1}{A(t)} \frac{\mathrm{d}A(t)}{\mathrm{d}t} \right)
  - \frac{1}{\epsilon} (x-X(t))^{\top} \left[
    \frac{1}{2} \frac{\mathrm{d}M(t)}{\mathrm{d}t}
    + M(t) \nabla_x \alpha_{i_0}(t,X(t))
  \right] (x-X(t)) = \re \gamma_{i_0 i_0} (t,X(t)), \\
\begin{split}
& \im \left( \frac{1}{A(t)} \frac{\mathrm{d}A(t)}{\mathrm{d}t} \right)
  - \frac{1}{\epsilon} (x-X(t))^{\top} \left[
    \frac{1}{2} \frac{\mathrm{d}N(t)}{\mathrm{d}t} -
    \frac{1}{2} \nabla_x^2 (P(t) \cdot \alpha_{i_0})(t,X(t)) +
    N(t) \nabla_x \alpha_{i_0}(t,X(t))
  \right] (x-X(t)) \\
& \qquad = \frac{1}{2} (x-X(t))^{\top} \nabla_x^2 \beta_{i_0}(t,X(t))
   (x-X(t)) + \im \gamma_{i_0 i_0}(t,X(t)).
\end{split}
\end{align*}
It is not difficult to find that the above equations hold if the
parameters of $G^{\epsilon}$ evolve according to the following system
of ODEs:
\begin{equation} \label{eq:gb_odes}
\begin{aligned}
& \frac{\mathrm{d}X(t)}{\mathrm{d}t} = \alpha_{i_0}(t,X(t)), \\
& \frac{\mathrm{d}S(t)}{\mathrm{d}t} = \beta_{i_0}(t,X(t)), \\
& \frac{\mathrm{d}P(t)}{\mathrm{d}t} = \nabla_x \beta_{i_0}(t,X(t))
  - [\nabla_x \alpha_{i_0}(t,X(t))]^{\top} P(t), \\
& \frac{\mathrm{d}A(t)}{\mathrm{d}t} =
  \gamma_{i_0 i_0}(t,X(t)) A(t), \\
& \frac{\mathrm{d}M(t)}{\mathrm{d}t} =
  - M(t) \nabla_x \alpha_{i_0}(t,X(t))
  - [\nabla_x \alpha_{i_0}(t,X(t))]^{\top} M(t), \\
& \frac{\mathrm{d}N(t)}{\mathrm{d}t} =
  \nabla_x^2 (P(t) \cdot \alpha_{i_0})(t,X(t))
  - \nabla_x^2 \beta_{i_0}(t,X(t))
  - N(t) \nabla_x \alpha_{i_0}(t,X(t))
  - [\nabla_x \alpha_{i_0}(t,X(t))]^{\top} N(t).
\end{aligned}
\end{equation}
The following proposition ensures that the matrix $M(t)$ is positive
definite:
\begin{proposition}
Suppose $B_1(t)$ and $B_2(t)$ are maps from $\mathbb{R}^+$ to
$\mathbb{R}^{m \times m}$, and both matrix functions are bounded and
piecewise continuous. Then the solution of
\begin{equation}
\begin{aligned}
& \frac{\mathrm{d}M(t)}{\mathrm{d}t} =
  B_1(t) + M(t) B_2(t) + [B_2(t)]^{\top} M(t), \\
& M(0) = M_0 \in \mathbb{R}^{m\times m}
\end{aligned}
\end{equation}
exists and is symmetric and positive definite for all $t > 0$ if the
following two conditions hold:
\begin{itemize}
\item $M_0$ is symmetric and positive definite;
\item $B_1(t)$ is symmetric and positive semi-definite for all $t > 0$.
\end{itemize}
\end{proposition}
The proof of the above proposition can be found in \cite{Dieci1994}.
We have thus established the Gaussian beam method for \eqref{eq:gb}.
One only needs to solve the ordinary differential system
\eqref{eq:gb_odes} to get the approximation solution of \eqref{eq:gb}.
Since all the terms with orders higher than or equal to
$O(\epsilon^{1/2})$ are dropped in the derivation of the method, the
error of the Gaussian beam method is $O(\epsilon^{1/2})$, which can be
made rigorous and we omit the details here.

\subsection{Surface hopping method for the ODE system} \label{sec:sh}

Note that for each fixed $x$, equation \eqref{eq:sh} is a system of
ODEs, while such ODEs can be solved easily by taking the exponent of
the matrix coefficients, in this section, we  propose an
alternative stochastic algorithm for such equations. The reason is
that the stochastic algorithm (the surface hopping algorithm) can be
easily combined with the Gaussian beam, while it is impossible for the
direct ODE solution as seen below. 

Define the matrix functions
$\Gamma(t,x) = (\tilde{\gamma}_{ij}(t,x))_{n\times n}$ and
$\Theta(t,x) = (\Theta_{ij}(t,x))_{n\times n}$ as
\begin{equation} \label{eq:tilde_gamma}
\tilde{\gamma}_{ij}(t,x) = \begin{cases}
  \gamma_{ij}(t,x) & \text{if } i \neq j, \\
  0 & \text{if } i = j,
\end{cases} \qquad
\Theta(t,x) = \exp \left( \int_0^t \Gamma(s,x) \mathrm{d}s \right),
\end{equation}
where the function $\exp(\cdot)$ means the matrix exponent. It is easy
to see that the solution to \eqref{eq:sh} is 
\begin{equation}
u_i^{\epsilon}(t,x) = \Theta_{ii_0}(t,x) v_{i_0}^{\epsilon}(x) =
  \Theta_{ii_0}(t,x) G^{\epsilon}(M_0, N_0, X_0, W_0, S_0, A_0; x),
\qquad i=1,\cdots,n.
\end{equation}
This solution does not have a Gaussian shape in general due to the
dependence of $\Theta_{ii_0}$ on $x$. To move forward, as we aim to
combine this with the Gaussian beam method, we would like to find an
approximation to the above solution with a Gaussian beam form. This
can be easily achieved by defining
\begin{equation}
U_i^{\epsilon}(t,x) = \Theta_{ii_0}(t,X_0) v_{i_0}^{\epsilon}(x) =
  \Theta_{ii_0}(t,X_0) G^{\epsilon}(M_0, N_0, X_0, W_0, S_0, A_0; x),
\qquad i=1,\cdots,n.
\end{equation}
Note that the argument of $\Theta$ has been replaced to $X_0$ above,
which is the fixed center of the Gaussian.  Using \eqref{eq:order}, we
have that
$U_i^{\epsilon}(t,x) - u_i^{\epsilon}(t,x) \sim O(\epsilon^{1/2})$.
Thus $U_i^{\epsilon}$ approximates the solution of \eqref{eq:sh} with
accuracy $O(\epsilon^{1/2})$, which accords with the accuracy for the
scalar equations in the last section. It is not difficult to find that
$U_i^{\epsilon}$ is the solution of the following linear system:
\begin{subequations} \label{eq:sh_U}
\begin{align}
\label{eq:ode_U}
& \frac{\partial U_i^{\epsilon}(t,x)}{\partial t} =
  \sum_{\substack{j=1\\ j\neq i}}^n
    \gamma_{ij}(t,X_0) U_j^{\epsilon}(t,x),
\quad i=1,\cdots,n, \\
& U_i^{\epsilon}(0,x) = v_i^{\epsilon}(x) =
  G^{\epsilon}(M_0, N_0, X_0, W_0, S_0, A_0; x) \delta_{ii_0},
\quad i=1,\cdots,n.
\end{align}
\end{subequations}

While $U_i^{\epsilon}$ gives an approximate solution, the solution
consists of a Gaussian for each surface $i = 1, \cdots, n$.  This is
undesirable, as if we combine this with the Gaussian beam method, each
time step would then create multiple Gaussian beams for each Gaussian,
which leads to exponential increase of the number of Gaussians as time
proceeds, which is impractical as stated in the beginning of this
section. To overcome this difficulty, below we
will propose a stochastic representation to the solution of
\eqref{eq:sh_U}, which is the core of the surface hopping method.  For
conciseness of notation, we abbreviate $\gamma_{ij}(t,X_0)$ as
$\gamma_{ij}(t)$ in the remainder of this section.

For the surface hopping method, we define a stochastic jump process
\begin{equation} \label{eq:ell}
\ell = \{ l_t : t \in \mathbb{R}^+ \}
\end{equation}
in the state space $\{1, \cdots, n\}$. We assume that the initial
state is fixed as $l_0 = i_0$. If the solution of the equation can be
formulated as
\begin{equation} \label{eq:exp}
U_i^{\epsilon}(t,x) = \mathbb{E}[\mathcal{F}_i(\ell; t,x)], \qquad
  i=1,\cdots,n
\end{equation}
for some functionals $\mathcal{F}_i(\cdot; \cdot, \cdot)$, then
$U_i^{\epsilon}$ can be numerically obtained by Monte Carlo sampling
$\ell$ and taking the average of $\mathcal{F}_i(\ell; t,x)$. Before we
determine the functionals $\mathcal{F}_i$, we would first like to give
a clear definition of the jump process $\ell$, and write out an
explicit formula for the expectation $U_i^{\epsilon}(t,x)$.

Assume the process $\ell$ has state $j$ at time $t$, and as $h
\rightarrow 0$, we have
\begin{equation} \label{eq:inf_prob}
\mathbb{P}(l_{t+h} = i \mid l_t = j) =
  \delta_{ij} + \Omega_{ij}(t) h + o(h),
\qquad i,j = 1,\cdots,n,
\end{equation}
where the infinitesimal hopping probability $\Omega_{ij}(t)$ satisfies
\begin{equation} \label{eq:Omega_ij}
\Omega_{ij}(t) \geqslant 0 \quad \text{if } i \neq j
\qquad \text{and} \qquad
\Omega_{jj}(t) = -\sum_{\substack{i=1\\i\neq j}}^n \Omega_{ij}(t),
\end{equation}
so that
\begin{equation}
\sum_{i = 1}^n \mathbb{P}(l_{t+h} = i \mid l_t = j) = 1.
\end{equation}
Let $J_{t_1, t_2}$ be the number of jumps in the time interval $(t_1,
t_2)$ (so it is a random variable), and then $\mathbb{P} (J_{t_1,t_2} = 0)$ is the probability that
no jump occurs in $(t_1,t_2)$, which can be obtained from
\eqref{eq:inf_prob} as
\begin{equation} \label{eq:J}
\mathbb{P} (J_{t_1,t_2} = 0) = \lim_{h\rightarrow 0}
  \prod_{k=1}^{\lceil (t_2-t_1)/h \rceil}
  \mathbb{P} (l_{t_1+kh}=l_{t_1} \mid l_{t_1+(k-1)h}=l_{t_1})
= \exp \left(
  \int_{t_1}^{t_2} \Omega_{l_{t_1}l_{t_1}}(s) \,\mathrm{d}s
\right).
\end{equation}
Based on this, we can show that
\begin{equation} \label{eq:exp_series}
\begin{split}
\mathbb{E}[\mathcal{F}_i(\ell;t,x)] &= \sum_{K=0}^{+\infty}
  \sum_{\substack{i_1=1\\i_1\neq i_0}}^n \cdots
  \sum_{\substack{i_K=1\\i_K\neq i_{K-1}}}^n \\
& \quad \times \int_0^t \!\! \int_0^{t_K} \!\! \cdots \! \int_0^{t_2}
    \mathcal{F}_i(\ell;t,x)
    \left( \prod_{k=1}^K \Omega_{i_k i_{k-1}}(t_k) \right)
    \exp \left( \int_0^t \Omega_{l_s l_s}(s) \,\mathrm{d}s \right)
  \,\mathrm{d}t_1 \cdots \,\mathrm{d}t_{K-1} \,\mathrm{d}t_K,
\end{split}
\end{equation}
where $\ell$ is a corresponding realization of the jump process
\eqref{eq:ell} given  $i_1, \cdots, i_K$ and $t_1, \cdots, t_K$ as
follows:
\begin{equation} \label{eq:jump}
l_s = i_k \quad
  \text{if} \quad s\in (t_k,t_{k+1}), \qquad k = 0,\cdots,K.
\end{equation}
In \eqref{eq:jump}, we take the convention that $t_0 = 0$ and $t_{K+1}
=t$ for simplicity. The summations and integrals \eqref{eq:exp_series}
shows that all possiblities of the jump process $\ell$ are taken into
account when evaluating the expectation of $\mathcal{F}_i$. The
product term contributes to the $K$ jumps, and the exponential term
contributes to the ``no-jump'' part of the process as in \eqref{eq:J}.
A detailed derivation of \eqref{eq:exp_series} can be found in the
Appendix \ref{sec:appendix}.

Next, we just need to write the solution of the system \eqref{eq:sh_U}
in the form of \eqref{eq:exp_series}. To do this, we write
$U_i^{\epsilon}$ as
\begin{equation} \label{eq:U_i}
U_i^{\epsilon}(t,x) = G^{\epsilon}(0;x) \delta_{ii_0} +
  \sum_{\substack{j_1=1\\ j_1 \neq i}}^n
    \int_0^t \gamma_{i j_1}(t_1) U_{j_1}^{\epsilon}(t_1, x)
    \,\mathrm{d}t_1.
\end{equation}
This equation is obtained simply by doing a time integration on both
sides of \eqref{eq:ode_U} over $[0,t]$, and it gives the integral form
of \eqref{eq:sh_U}. In order to get a series, we insert \eqref{eq:U_i}
into the right hand side of \eqref{eq:U_i} itself, and get a new
expression of $U_i^{\epsilon}$: 
\begin{equation}
\begin{split}
U_i^{\epsilon}(t,x) &= G^{\epsilon}(0;x) \delta_{ii_0} +
  \sum_{\substack{j_1=1\\ j_1 \neq i}}^n \int_0^t
    \gamma_{ij_1}(t_1) G^{\epsilon}(0;x) \delta_{j_1 i_0}
  \,\mathrm{d}t_1 + \sum_{\substack{j_1=1\\ j_1 \neq i}}^n
    \sum_{\substack{j_2=1\\ j_2 \neq j_1}}^n \int_0^t \! \int_0^{t_1}
      \gamma_{i j_1}(t_1) \gamma_{j_1 j_2}(t_2)
      U_{j_2}^{\epsilon}(t_1, x)
    \,\mathrm{d}t_2 \,\mathrm{d}t_1 \\
&= G^{\epsilon}(0;x) \delta_{ii_0} +
  \sum_{\substack{i_1=1\\ i_1 \neq i_0}}^n \int_0^t
    \gamma_{i_1 i_0}(t_1) G^{\epsilon}(0;x) \delta_{ii_1}
  \,\mathrm{d}t_1 + \sum_{\substack{j_1=1\\ j_1 \neq i}}^n
    \sum_{\substack{j_2=1\\ j_2 \neq j_1}}^n \int_0^t \! \int_0^{t_1}
      \gamma_{i j_1}(t_1) \gamma_{j_1 j_2}(t_2)
      U_{j_2}^{\epsilon}(t_1, x)
    \,\mathrm{d}t_2 \,\mathrm{d}t_1.
\end{split}
\end{equation}
By repeatedly inserting \eqref{eq:U_i} into the newly obtained
expression of $U_i^{\epsilon}$ and applying the following change of
variables
\begin{equation}
\begin{split}
& \sum_{\substack{j_1 = 1\\ j_1 \neq i}}^n
\sum_{\substack{j_2 = 1\\ j_2 \neq j_1}}^n \cdots
\sum_{\substack{j_K = 1\\ j_K \neq j_{K-1}}}^n
  \int_0^t \int_0^{t_1} \cdots \int_0^{t_{K-1}}
    \gamma_{ij_1}(t_1) \gamma_{j_1 j_2}(t_2) \cdots
    \gamma_{j_{K-1} j_K}(t_K) G^{\epsilon}(0;x) \delta_{j_K i_0}
  \,\mathrm{d}t_K \cdots \,\mathrm{d}t_2 \,\mathrm{d}t_1 \\
={} & \sum_{\substack{i_1 = 1\\ i_1 \neq i_0}}^n
\sum_{\substack{i_2 = 1\\ i_2 \neq i_1}}^n \cdots \!\!
\sum_{\substack{i_K = 1\\ i_K \neq i_{K-1}}}^n
  \int_0^t \!\! \int_0^{t_K} \cdots \! \int_0^{t_2}
    \gamma_{i_K i_{K-1}}(t_K) \gamma_{i_{K-1} i_{K-2}}(t_{K-1}) \cdots
    \gamma_{i_1 i_0}(t_1) G^{\epsilon}(0;x) \delta_{ii_K}
  \,\mathrm{d}t_1 \cdots \,\mathrm{d}t_{K-1} \,\mathrm{d}t_K,
\end{split}
\end{equation}
we will eventually obtain the following series form of
$U_i^{\epsilon}$:
\begin{equation} \label{eq:Ui_series}
U_i^{\epsilon}(t,x) = \sum_{K=0}^{+\infty}
  \sum_{\substack{i_1 = 1\\ i_1 \neq i_0}}^n
  \sum_{\substack{i_2 = 1\\ i_2 \neq i_1}}^n \cdots
  \sum_{\substack{i_K = 1\\ i_K \neq i_{K-1}}}^n
    \int_0^t \int_0^{t_K} \cdots \int_0^{t_2}
      \left( \prod_{k=1}^K \gamma_{i_k i_{k-1}}(t_k) \right)
      G^{\epsilon}(0;x) \delta_{ii_K}
    \,\mathrm{d}t_1 \cdots \,\mathrm{d}t_{K-1} \,\mathrm{d}t_K.
\end{equation}
Now comparing \eqref{eq:Ui_series} with \eqref{eq:exp_series}, one
easily finds that \eqref{eq:exp} is fulfilled by setting
\begin{equation} \label{eq:F}
\mathcal{F}_i(\ell;t,x) = \exp \left(
  -\int_0^t \Omega_{l_s l_s}(s) \,\mathrm{d}s
\right) \prod_{k=1}^K
  \frac{\gamma_{i_k i_{k-1}}(t_k)}{\Omega_{i_k i_{k-1}}(t_k)}
  G^{\epsilon}(0;x) \delta_{ii_K}.
\end{equation}
Note that when $i = j$, $\Omega_{ij}$ is defined in \eqref{eq:Omega_ij}, and
therefore it remains only to choose the infinitesimal jump probabilities
$\Omega_{ij}(t)$ for $i \neq j$. The choice is nonunique. In fact, any positive
function independent of $\epsilon$ can be used. In our algorithm, we use
\begin{equation} \label{eq:jump_prob}
\Omega_{ij}(t) = |\gamma_{ij}(t)|, \qquad
  i,j = 1,\cdots,n \quad \text{and} \quad i \neq j,
\end{equation}
so that the hop does not change the magnitude of the amplitude of the Gaussian
beam. Based on \eqref{eq:F}, we can define a function
$\omega(t)$ as the solution of
\begin{equation} \label{eq:ode_Omega}
\frac{\mathrm{d}\omega(t)}{\mathrm{d}t} = -\Omega_{l_t l_t}(t) +
  \mathrm{i} \frac{\mathrm{d} J_{0,t}}{\mathrm{d} t}
    \arg \gamma_{l_{t+} l_{t-}}(t) =
\sum_{\substack{j=1\\j\neq l_t}}^m |\gamma_{jl_t}(t)| +
  \mathrm{i} \frac{\mathrm{d} J_{0,t}}{\mathrm{d} t}
    \arg \gamma_{l_{t+} l_{t-}}(t), \qquad \omega(0) = 0,
\end{equation}
where  $l_{t-}$ ($l_{t+}$) denotes the left (right) limit of $l_t$ at
time $t$. Since $J_{0,t}$ is the number of jumps before time $t$, it is right
continuous with left limits, and for a given jump process defined by
\eqref{eq:jump}, the derivative of $J_{0,t}$ is (in the distributional sense)
\begin{equation}
\frac{\mathrm{d} J_{0,t}}{\mathrm{d}t} = \sum_{k=1}^K \delta(t - t_k).
\end{equation}
Therefore the solution of \eqref{eq:ode_Omega} is
\begin{equation}
\omega(t) = \exp \Bigg(
  \sum_{\substack{j=1\\j\neq l_t}}^m |\gamma_{jl_t}(t)|
\Bigg) \prod_{k=1}^K
  \frac{\gamma_{i_k i_{k-1}}(t_k)}{|\gamma_{i_k i_{k-1}}(t_k)|},
\end{equation}
and then we have
\begin{equation} \label{eq:sh_final}
U_i^{\epsilon}(t,x) =
  \mathbb{E}[\exp(\omega(t)) G^{\epsilon}(0;x) \delta_{il_t}].
\end{equation}
By now, it is clear that the numerical method for the system
\eqref{eq:sh_U} contains the following ingredients:
\vspace{-\topsep}
\begin{itemize}
\item Draw a number of realizations of the jump process $\ell$
  according to the infinitesimal jump probability \eqref{eq:Omega_ij}.
\item For each trajectory, solve $\omega(t)$ along the trajectory from
  \eqref{eq:ode_Omega}.
\item Take the average of all trajectories as \eqref{eq:sh_final} to
  get $U_i^{\epsilon}$.
\end{itemize}

Due to the ``$\delta_{il_t}$'' in the final result
\eqref{eq:sh_final}, we see that at time $t$, the jump process $l_t$
only contributes to the $l_t$th component of the solution. Therefore,
this algorithm can be described as a Gaussian hopping between
different surfaces while the jump in $l_t$ occurs. When there is no
hop and the Gaussian stays at the $i$th surface, the Gaussian just
evolves by changing its amplitude gradually (the coefficient
$\exp(\omega(t))$ in \eqref{eq:sh_final} is continuous when $l_t$ is a
constant), which looks like travelling along a predefined surface. At
the same time, the Gaussian is counted into the $i$th component of the
solution. However, when the Gaussian hops from the $i$th surface to
the $j$th surface, it stops its contribution to the $i$th component
and becomes a part of the solution of the $j$th component. Exactly at
the time point of such a hop, the amplitude of the Gaussian has a
sudden change due to the discontinuity in $J_{0,t}$. This is why the algorithm is
called surface hopping method. As we will see in the next subsection, when combined with the Gaussian beam method, the surface hopping provides a way to solve the high dimensional transport system. 

\subsection{The surface hopping Gaussian beam method} \label{sec:shgb}
We now combine the two methods introduced in the last two sections to
construct a stochastic solver for the original system
\eqref{eq:equations}. A natural way of such combination is to use the
algorithm described in Section \ref{sec:sh} to define a jump process,
and do the same to the Gaussian as in the surface hopping method when
a hop occurs. When there is no hop, the Gaussian $G^{\epsilon}$
evolves as in the Gaussian beam method described in Section
\ref{sec:gb}, and its amplitude changes as in Section \ref{sec:sh}. A
simple way to interpret such an combination is to use the idea of the
time splitting method, which solve equations \eqref{eq:gb} and
\eqref{eq:sh} alternately. Note that in the beginning of Section
\ref{sec:sh}, we freeze the coefficients at point $X_0$ since $X_0$ is
the center of the Gaussian (see \eqref{eq:ode_U}). Here, the Gaussian
beam is moving around, and therefore when solving \eqref{eq:sh}, it is
natural to freeze the coefficient $|\gamma_{ij}(t,x)|$ at the center
of the current Gaussian beam. At each time step, we need to generate a
random number to determine whether the trajectory hops at the current
time step. Below, we are going to carry out such an idea without
explicitly using the time splitting, so that much less random numbers
are needed in the implementation.

We define the jump process $l_t$ as \eqref{eq:inf_prob}, while the
infinitesimal jump probability \eqref{eq:jump_prob} is modified as
\begin{equation} \label{eq:shgb_omega}
\Omega_{ij}(t) = |\gamma_{ij}(t,X(t))|, \qquad
  i,j = 1,\cdots,n \quad \text{and} \quad i \neq j,
\end{equation}
where $X(t)$ is the center of the Gaussian as in Section \ref{sec:gb}.
For each $l_t$, we solve the following ordinary differential system:
\begin{equation} \label{eq:shgb_odes}
\begin{aligned}
& \frac{\mathrm{d}X(t)}{\mathrm{d}t} = \alpha_{l_t}(t,X(t)), \\
& \frac{\mathrm{d}S(t)}{\mathrm{d}t} = \beta_{l_t}(t,X(t)), \\
& \frac{\mathrm{d}P(t)}{\mathrm{d}t} = \nabla_x \beta_{l_t}(t,X(t))
  - [\nabla_x \alpha_{l_t}(t,X(t))]^{\top} P(t), \\
& \frac{\mathrm{d}A(t)}{\mathrm{d}t} =
  \gamma_{l_t l_t}(t,X(t)) A(t), \\
& \frac{\mathrm{d}M(t)}{\mathrm{d}t} =
  - M(t) \nabla_x \alpha_{l_t}(t,X(t))
  - [\nabla_x \alpha_{l_t}(t,X(t))]^{\top} M(t), \\
& \frac{\mathrm{d}N(t)}{\mathrm{d}t} =
  \nabla_x^2 (P(t) \cdot \alpha_{l_t})(t,X(t))
  - \nabla_x^2 \beta_{l_t}(t,X(t))
  - N(t) \nabla_x \alpha_{l_t}(t,X(t))
  - [\nabla_x \alpha_{l_t}(t,X(t))]^{\top} N(t), \\
& \frac{\mathrm{d}\omega(t)}{\mathrm{d}t} =
  \sum_{\substack{j=1\\j\neq l_t}}^m |\gamma_{jl_t}(t,X(t))| +
    \mathrm{i} \frac{\mathrm{d} J_{0,t}}{\mathrm{d} t}
      \arg \gamma_{l_{t+} l_{t-}}(t,X(t)),
\end{aligned}
\end{equation}
with initial conditions
\begin{equation} \label{eq:shgb_init}
X(0) = X_0, \quad S(0) = S_0, \quad P(0) = W_0, \quad
A(0) = A_0, \quad M(0) = M_0, \quad N(0) = N_0, \quad \omega(0) = 0.
\end{equation}
It is clear that the first six equations in \eqref{eq:shgb_odes} come
from the Gaussian beam method, but the index $i_0$ in
\eqref{eq:gb_odes} is replaced by $l_t$ in the above system, since the
Gaussian should follow the ``current'' surface which has the index
$l_t$. The last equation in \eqref{eq:shgb_odes} comes from the
surface hopping method for the ordinary differential system. Inspired
by \eqref{eq:sh_final}, an approximate solution of
\eqref{eq:equations} is expected to be
\begin{equation} \label{eq:shgb}
U_i^{\epsilon}(t,x) =
  \mathbb{E}[\exp(\omega(t)) G^{\epsilon}(t;x) \delta_{il_t}],
\qquad i=1,\cdots,n.
\end{equation}
Below we will verify that \eqref{eq:shgb} satisfies
\eqref{eq:equations} up to a residual $O(\epsilon^{1/2})$.

The derivation in Section \ref{sec:sh} shows that \eqref{eq:shgb} can
be expanded as \eqref{eq:Ui_series} with $G^{\epsilon}(0;x)$ replaced
by $G^{\epsilon}(t;x)$. For simplicity, we can employ the definition
\eqref{eq:tilde_gamma} to replace $\gamma_{i_k i_{k-1}}$ by
$\tilde{\gamma}_{i_k i_{k-1}}$ in \eqref{eq:Ui_series} and remove the
constraints $i_1 \neq i_0$, $i_2 \neq i_1$, $\cdots$, $i_K \neq
i_{K-1}$ in the summations. Furthermore, we let $\tau_K = (t_1, t_2,
\cdots, t_K)^{\top}$ and define
\begin{equation}
S_t^K = \{ \tau_K \mid
  0\leqslant t_1 \leqslant \cdots \leqslant t_K \leqslant t \}.
\end{equation}
Thus for any function $f(t,t_1,\cdots,t_K)$, one has
\begin{equation}
\int_{S_t^K} f(t, \tau_K) \,\mathrm{d}\tau_K =
  \int_0^t \int_0^{t_K} \cdots \int_0^{t_2} f(t, t_1, \cdots, t_K)
    \,\mathrm{d}t_1 \cdots \,\mathrm{d}t_{K-1} \,\mathrm{d}t_K.
\end{equation}
Using these notations, we can rewrite $U_i^{\epsilon}$ in a much simpler way as
\begin{equation} \label{eq:Ui_simplified}
U_i^{\epsilon}(t,x) = \sum_{K=0}^{+\infty} \sum_{i_1,\cdots,i_K=1}^n
  \int_{S_t^K} \left(
    \prod_{k=1}^K \tilde{\gamma}_{i_k i_{k-1}}(t_k, X(t_k))
  \right) G^{\epsilon}(t;x) \delta_{i i_K} \,\mathrm{d}\tau_K.
\end{equation}
Now we are going to take the time derivative of $U_i^{\epsilon}$.
Note that the first six equations in \eqref{eq:shgb_odes} form a closed
subsystem, and in the right hand side of the subsystem, discontinuities exist
only in the jump process. Therefore $G^{\epsilon}(t,x)$ is a continuous
function, and for a fixed jump process $l_t$, we can use the result in Section
\ref{sec:gb} to get
\begin{equation} \label{eq:G_eq}
\frac{\partial G^{\epsilon}(t;x)}{\partial t} +
  \alpha_{l_t}(t,x)^{\top} \nabla_x G^{\epsilon}(t;x) =
\left(
  \frac{\mathrm{i}}{\epsilon} \beta_{l_t}(t,x) + \gamma_{l_t l_t}(t,x)
\right) G^{\epsilon}(t;x) + O(\epsilon^{\frac{1}{2}}).
\end{equation}
Using
\begin{equation}
\frac{\mathrm{d}}{\mathrm{d}t} \int_{S_t^K} f(t,\tau_K)
\,\mathrm{d}\tau_K = \int_{S_t^{K-1}} f(t,\tau_{K-1},t)
\,\mathrm{d}\tau_{K-1} + \int_{S_t^K} \partial_t f(t,\tau_K)
\,\mathrm{d}\tau_K,
\end{equation}
we can take the time derivative of \eqref{eq:Ui_simplified} and get
\begin{equation}
\begin{split}
& \frac{\partial U_i^{\epsilon}(t,x)}{\partial t} =
  \sum_{K=0}^{+\infty} \sum_{i_1, \cdots, i_K=1}^n
  \tilde{\gamma}_{i_K i_{K-1}}(t,X(t)) \int_{S_t^{K-1}} \left(
    \prod_{k=1}^{K-1} \tilde{\gamma}_{i_k i_{k-1}}(t_k, X(t_k))
  \right) G^{\epsilon}(t;x) \delta_{ii_K} \,\mathrm{d}\tau_{K-1} + {} \\
& + \sum_{K=0}^{+\infty} \sum_{i_1, \cdots, i_K=1}^n \int_{S_t^K}
    \left( \prod_{k=1}^K \tilde{\gamma}_{i_k i_{k-1}}(t_k, X(t_k)) \right)
    \left[
      \left(
        \frac{\mathrm{i}}{\epsilon} \beta_{i_K}(t,x) + \gamma_{i_K i_K}(t,x)
        - \alpha_{i_K}(t,x)^{\top} \nabla_x
      \right) G^{\epsilon}(t;x)
    \right] \delta_{ii_K} \,\mathrm{d}\tau_K + O(\epsilon^{\frac{1}{2}}).
\end{split}
\end{equation}
This equation shows that in \eqref{eq:Ui_simplified}, the derivative
of the integrand provides the contribution on the current surface,
while the derivative of the integral bounds provides the contribution
from the surface hopping term. In detail, we can write the above
equality as $\partial_t U_i^{\epsilon} = \mathrm{I} + \mathrm{II} +
O(\epsilon^{\frac{1}{2}})$ and get
\begin{equation} \label{eq:I}
\begin{split}
\mathrm{I} &= \sum_{K=0}^{+\infty} \sum_{i_1, \cdots, i_{K-1}=1}^n
  \tilde{\gamma}_{i i_{K-1}}(t,X(t)) \int_{S_t^{K-1}} \left(
    \prod_{k=1}^{K-1} \tilde{\gamma}_{i_k i_{k-1}}(t_k,X(t_k))
  \right) G^{\epsilon}(t;x) \,\mathrm{d} \tau_{K-1} \\
& = \sum_{K=0}^{+\infty} \sum_{i_1, \cdots, i_{K-1}=1}^n
  \tilde{\gamma}_{i i_{K-1}}(t,x) \int_{S_t^{K-1}} \left(
    \prod_{k=1}^{K-1} \tilde{\gamma}_{i_k i_{k-1}}(t_k,X(t_k))
  \right) G^{\epsilon}(t;x) \,\mathrm{d} \tau_{K-1}
    + O(\epsilon^{\frac{1}{2}}) \\
&= \sum_{j=1}^n \sum_{K=0}^{+\infty} \sum_{i_1, \cdots, i_{K-1}=1}^n
  \tilde{\gamma}_{ij}(t,x) \int_{S_t^{K-1}} \left(
    \prod_{k=1}^{K-1} \tilde{\gamma}_{i_k i_{k-1}}(t_k, X(t_k))
  \right) G^{\epsilon}(t;x) \delta_{ji_{K-1}} \,\mathrm{d} \tau_{K-1}
    + O(\epsilon^{\frac{1}{2}}) \\
&= \sum_{j=1}^n \tilde{\gamma}_{ij}(t,x) U_j^{\epsilon}(t,x)
  + O(\epsilon^{\frac{1}{2}}),
\end{split}
\end{equation}
where \eqref{eq:order} is used in the second equality. For
$\mathrm{II}$, it is obvious that
\begin{equation} \label{eq:II}
\mathrm{II} = \left(
  \frac{\mathrm{i}}{\epsilon} \beta_i(t,x) +
    \gamma_{ii}(t,x) - \alpha_i(t,x)^{\top} \nabla_x
\right) U_i^{\epsilon}(t,x).
\end{equation}
Adding \eqref{eq:I} and \eqref{eq:II}, we conclude that
$U_i^{\epsilon}(t,x)$ defined by \eqref{eq:shgb} satisfies
\begin{equation}
\frac{\partial U_i^{\epsilon}(t,x)}{\partial t} +
  \alpha_i(t,x)^{\top} \nabla_x U_i^{\epsilon}(t,x) =
\frac{\mathrm{i}}{\epsilon} \beta_i(t,x) U_i^{\epsilon}(t,x) +
  \sum_{j=1}^n \gamma_{ij}(t,x) U_i^{\epsilon}(t,x) +
  O(\epsilon^{\frac{1}{2}}),
\end{equation}
and therefore is an approximate solution of \eqref{eq:equations}.

The algorithm to obtain such an approximate solution $U_i^{\epsilon}$
is implied in \eqref{eq:shgb_omega}--\eqref{eq:shgb}, and it will be
made more explicit in the next section. The equations
\eqref{eq:shgb_odes} show that while the Gaussian beam is traveling
along the surfaces, it changes not only the amplitude, but also its
shape, center and oscillation frequency. Thus, most of the time, the
Gaussian beam moves following a ``trajectory'' on a surface, and at
several time points, it hops from one surface to another
stochastically. The average of all the trajectories gives the solution
of the system \eqref{eq:equations}.

\subsection{Outline of the algorithm} \label{sec:outline}
Based on the above derivation, we list the detailed steps of the
algorithm in this section. Before that, we claim that the time of the
first jump $t_1$ can be picked by generating a random number $Y \sim
U(0,1)$, and solve
\begin{equation} \label{eq:hop_time}
Y = 1 -
  \exp \left( \int_0^{t_1} \Omega_{i_0 i_0}(s) \,\mathrm{d}s \right)
= 1 - \exp \left( -\int_0^{t_1}
  \sum_{k=1}^{n} |\tilde{\gamma}_{k i_0}(s, X(s))|
\,\mathrm{d}s \right).
\end{equation}
This can be observed from \eqref{eq:J}, which tells us that the
probability $\mathbb{P}(t_1 < T)$ is $1 - \mathbb{P}(J_{0,T} = 0)$,
and therefore $Y \sim U(0,1)$ if $Y = 1 - \mathbb{P}(J_{0,t_1} = 0)$.
Since the right hand side of \eqref{eq:hop_time} is an increasing
function of $t_1$, the solution is unique, and the equation can be
solved while evolving the trajectory. Other jump times $t_2$, $t_3$,
$\cdots$ can be similarly generated.

Now we are ready to state our algorithm. We assume that we want to
obtain $U_i^{\epsilon}(T,x)$ for all $j = 1,\cdots,n$.
\begin{enumerate}
\item Set $N_{\mathrm{traj}} \leftarrow 0$ and $U_j^{\epsilon}(T,x)
\leftarrow 0$ for all $j = 1,\cdots,n$.
\item \label{item:traj_init} Set $i \leftarrow i_0$, $t \leftarrow 0$,
$\tilde{\omega} \leftarrow 0$, $\mathit{flag} \leftarrow
\mathit{false}$. Set the initial conditions \eqref{eq:shgb_init}.
Generate a random number $Y \sim U(0,1)$.
\item \label{item:evolve} Select the time step $\Delta t$. Assume that
no hops occurs in $(t, t + \Delta t)$, and solve $\omega(t+\Delta t)$
according to \eqref{eq:shgb_odes}.
\item \label{item:hop_time} If $\exp(\tilde{\omega} - \omega(t+\Delta
t)) \leqslant 1 - Y$, then solve the equation of the jump time
$\tilde{t}$
\begin{displaymath}
\exp(\tilde{\omega} - \omega(\tilde{t})) = 1 - Y
\end{displaymath}
by interpolation of $\omega(t)$ in $[t, t+\Delta t]$, and set $\Delta
t \leftarrow \tilde{t} - t$ and $\mathit{flag} \leftarrow
\mathit{true}$.
\item \label{item:hop} Solve the whole system \eqref{eq:shgb_odes} by
$\Delta t$ with the assumption that no hops occurs in $(t, t + \Delta
t)$. Set $t \leftarrow t + \Delta t$. If $\mathit{flag}$ is
$\mathit{true}$, then generate a random integer $j \in \{1, \cdots,
n\} \backslash \{i\}$ from the probability mass function
\begin{displaymath}
f(j) = |\gamma_{ji}(t,X(t))| \bigg/
  \sum_{\substack{k=1\\k\neq i}}^n |\gamma_{ki}(t,X(t))|,
\qquad j \in \{1,\cdots,n\} \backslash \{i\},
\end{displaymath}
generate a new $Y \sim U(0,1)$, and set
\begin{displaymath}
\omega(t) \leftarrow \omega(t) + \mathrm{i} \arg \gamma_{ji}(t,X(t)),
\qquad i \leftarrow j, \qquad \tilde{\omega} \leftarrow \omega(t),
\qquad \mathit{flag} \leftarrow \mathit{false}.
\end{displaymath}
\item If $t < T$, return to step \ref{item:evolve}. If $t \geqslant
T$, then set
\begin{displaymath}
U_i^{\epsilon}(T,x) \leftarrow U_i^{\epsilon}(T,x) +
  \exp(\omega(T)) G^{\epsilon}(T;x), \qquad
N_{\mathrm{traj}} \leftarrow N_{\mathrm{traj}} + 1.
\end{displaymath}
\item If $N_{\mathrm{traj}}$ is large enough, then set
$U_j^{\epsilon}(T,x) \leftarrow U_j^{\epsilon}(T,x) /
N_{\mathrm{traj}}$ for all $j = 1,\cdots,n$ and stop. Otherwise,
return to step \ref{item:traj_init}.
\end{enumerate}
In this algorithm, the jump process $l_t$ is hidden in the parameter
$i$. The boolean $\mathit{flag}$ indicates whether the hopping occurs
at the current time step. The parameter $\tilde{\omega}$ records the
value of $\omega$ at the time of last jump, and thus according to
\eqref{eq:hop_time}, the equation in step \ref{item:hop_time}
determines the time of the first jump since the last jump. When
it is known that a hop occurs between $t$ and $t + \Delta t$, we
shrink the time step $\Delta t$ such that $t + \Delta t$ is exactly
the time point of the next jump. The surface hopping is applied at step
\ref{item:hop}. It is obvious that we don't need to generate a random
number at every time step as in a naive time splitting method. In our
implementation, the system \eqref{eq:shgb_odes} is solved with a
third-order Runge-Kutta method in step \ref{item:evolve} and
\ref{item:hop}.

\begin{remark*}
The above method can also be used to solve the following extension of
the system \eqref{eq:equations}:
\begin{equation} \label{eq:extended}
\frac{\partial u_i^{\epsilon}(t,x)}{\partial t} +
  \alpha_i(t,x)^{\top} \nabla_x u_i^{\epsilon}(t,x) =
\frac{\mathrm{i}}{\epsilon} \beta_i(t,x) u_i^{\epsilon}(t,x) +
  \sum_{j=1}^n \gamma_{ij}(t,x) u_j^{\epsilon}(t,x) +
  \sum_{j=1}^n \nu_{ij}(t,x) \overline{u_j^{\epsilon}(t,x)},
\quad i=1,\cdots,n.
\end{equation}
To do this, one can define $w_i^{\epsilon}(t,x) =
\overline{u_i^{\epsilon}(t,x)}$ and rewrite the system
\eqref{eq:extended} as
\begin{displaymath}
\begin{aligned}
\frac{\partial u_i^{\epsilon}(t,x)}{\partial t} +
  \alpha_i(t,x)^{\top} \nabla_x u_i^{\epsilon}(t,x) & =
\frac{\mathrm{i}}{\epsilon} \beta_i(t,x) u_i^{\epsilon}(t,x) +
  \sum_{j=1}^n \gamma_{ij}(t,x) u_j^{\epsilon}(t,x) +
  \sum_{j=1}^n \nu_{ij}(t,x) w_j^{\epsilon}(t,x),
\quad i=1,\cdots,n, \\
\frac{\partial w_i^{\epsilon}(t,x)}{\partial t} +
  \alpha_i(t,x)^{\top} \nabla_x w_i^{\epsilon}(t,x) & =
\frac{\mathrm{i}}{\epsilon}
  \overline{\beta_i(t,x)} w_i^{\epsilon}(t,x) +
  \sum_{j=1}^n \overline{\gamma_{ij}(t,x)} w_j^{\epsilon}(t,x) +
  \sum_{j=1}^n \overline{\nu_{ij}(t,x)} u_j^{\epsilon}(t,x),
\quad i=1,\cdots,n.
\end{aligned}
\end{displaymath}
The above system again has the form of \eqref{eq:equations}, and thus
the algorithm can be applied. We can average the results of
$u_i^{\epsilon}(t,x)$ and $\overline{w_i^{\epsilon}(t,x)}$ to get the
final result. Due to the special relation between $u_i^{\epsilon}$ and
$w_i^{\epsilon}$, when implementing the algorithm, we do not need to
double the number of unknowns. When a hop from a
$u_i^{\epsilon}$-surface to a $w_j^{\epsilon}$-surface occurs, we just
need to set
\begin{equation} \label{eq:switch}
A(t) \leftarrow \overline{A(t)}, \quad
S(t) \leftarrow -S(t), \quad
P(t) \leftarrow -P(t), \quad
N(t) \leftarrow -N(t), \quad
\omega(t) \leftarrow \overline{\omega(t)}
\end{equation}
after the hop, and then continue the algorithm as if the Gaussian beam
is on the $u_j^{\epsilon}$-surface. Thus $w_i^{\epsilon}$ does not
appear in the algorithm, while \eqref{eq:switch} is a hidden switch
between $u$ and $w$.
\end{remark*}

\subsection{Algorithm for general initial data} \label{sec:init_data}
The above discussion is based on the initial condition
\eqref{eq:init}, which is a single Gaussian beam. This algorithm can
be naturally extended to more general case due to the linearity of the
equations. Suppose the initial data can be written as
\begin{equation} \label{eq:init_sample}
v_i^{\epsilon}(x) = \int_{\mathbb{R}^m} \int_{\mathbb{R}^m}
  \int_{\mathbb{R}^{m\times m}} \int_{\mathbb{R}^{m\times m}}
    f_i^{\epsilon}(M, N, X, P)
    G^{\epsilon}(M, N, X, P, S_i(M,N,X,P), A_i(M,N,X,P); x)
  \,\mathrm{d}M \,\mathrm{d}N \,\mathrm{d}X \,\mathrm{d}P
\end{equation}
for every $i = 1,\cdots,n$. If $f_i^{\epsilon}(M,N,X,P)$ has a support
with nonzero measure and is absolutely integrable, then we can select
$f_i^{\epsilon}(M,N,X,P)$ and $A_i(M,N,X,P)$ appropriately such that
$f_i^{\epsilon}(M,N,X,P) \geqslant 0$ and
\begin{equation} \label{eq:normal}
\int_{\mathbb{R}^m} \int_{\mathbb{R}^m}
  \int_{\mathbb{R}^{m\times m}} \int_{\mathbb{R}^{m\times m}}
    f_i^{\epsilon}(M, N, X, P)
  \,\mathrm{d}M \,\mathrm{d}N \,\mathrm{d}X \,\mathrm{d}P = 1.
\end{equation}
In detail, if the positivity of $f_i^{\epsilon}(M,N,X,P)$ or
\eqref{eq:normal} is not satisfied, we can set
\begin{equation}
f_i^{\epsilon}(M,N,X,P) \leftarrow |f_i^{\epsilon}(M,N,X,P)| / C_i,
\qquad A_i(M,N,X,P) \leftarrow C_i f_i^{\epsilon}(M,N,X,P) /
|f_i^{\epsilon}(M,N,X,P)|,
\end{equation}
where $C_i$ is selected such that \eqref{eq:normal} can be fulfilled.

To solve the system \eqref{eq:equations} with initial conditions
\eqref{eq:init_sample}, we only need to sample the initial parameters
of the Gaussian beam according to the density functions
$f_i^{\epsilon}$. Concretely speaking, let
\begin{equation}
\mathcal{I} =
  \{ i \mid \mathrm{supp} \, v_i^{\epsilon} \neq \varnothing \},
\end{equation}
and then we just need to change step \ref{item:traj_init} to the
following:
\begin{enumerate}
\setcounter{enumi}{1}
\item Randomly pick $i$ from $\mathcal{I}$. Draw a sample $M_0, N_0,
X_0, W_0$ from $f_i^{\epsilon}(M,N,X,P)$ and set
\begin{displaymath}
S_0 \leftarrow S_i(M_0, N_0, X_0, W_0), \quad
A_0 \leftarrow n_{\mathcal{I}} A_i(M_0, N_0, X_0, W_0), \quad
t \leftarrow 0, \quad \tilde{\omega} \leftarrow 0,
\quad \mathit{flag} \leftarrow \mathit{false}.
\end{displaymath}
Generate a random number $Y \sim U(0,1)$.
\end{enumerate}
Thus all the trajectories start from different surfaces and different
positions. Here $n_{\mathcal{I}}$ is the number of elements in the set
$\mathcal{I}$. It is multiplied to the initial amplitude $A_0$ since
when selecting the initial surface, the probablity of picking any
element from $\mathcal{I}$ is $1/n_{\mathcal{I}}$, which will be
introduced as a factor to the solution when we calculate the
expectation.  The expectation of these Gaussian beams is the
expectation of the system's approximate solutions with initial data
sampled in the above step, which is just an approximate solution of
the system \eqref{eq:equations} due to its linearity.

The initial value decomposition and its accuracy are also discussed in
a number of papers. Some recent development includes
\cite{Tanushev2008,Qian2010, Ariel2011}.

\section{Numerical examples} \label{sec:num_ex}
\subsection{An example without modelling error}
We first consider the following two-dimensional system:
\begin{equation} \label{eq:ex1}
\begin{aligned}
& \frac{\partial u_1^{\epsilon}}{\partial t} -
  (x_1 + x_2 + 1) \frac{\partial u_1^{\epsilon}}{\partial x_1} -
  (2x_1 - x_2 - 1) \frac{\partial u_1^{\epsilon}}{\partial x_2} =
\frac{\mathrm{i}}{\epsilon} |x|^2 u_1^{\epsilon}
  - u_1^{\epsilon} + \frac{1}{2} u_2^{\epsilon}, \\
& \frac{\partial u_2^{\epsilon}}{\partial t} -
  \left( \frac{1}{2} x_1 + x_2 - 1 \right)
    \frac{\partial u_2^{\epsilon}}{\partial x_1} -
  \left( x_1 - \frac{1}{2} x_2 + 1 \right)
    \frac{\partial u_1^{\epsilon}}{\partial x_2} =
\frac{\mathrm{i}}{\epsilon} (x_1 + x_2)^2 u_2^{\epsilon}
  + \frac{2}{3} u_1^{\epsilon} -2 u_2^{\epsilon},
\end{aligned}
\end{equation}
with initial value
\begin{equation}
u_1^{\epsilon}(0,x) = \exp \left( -\frac{|x|^2}{2} \right), \qquad
u_2^{\epsilon}(0,x) = 0.
\end{equation}
In order to apply the Gaussian beam method, we write the initial value
of $u_1^{\epsilon}$ as
\begin{equation}
u_1^{\epsilon}(0,x) = \int_{\mathbb{R}^2}
  f^{\epsilon}(X) \cdot \frac{1}{\epsilon} \exp \left(
    -\frac{|x-X|^2}{2\epsilon}
  \right) \,\mathrm{d}X,
\end{equation}
where $f^{\epsilon}(X)$ is a normal distribution
\begin{equation}
f^{\epsilon}(X) = \frac{1}{2\pi(1-\epsilon)} \exp
  \left( -\frac{|X|^2}{2(1-\epsilon)} \right).
\end{equation}
Thus the initial Gaussian beam can be obtained by drawing $X$ from
$f^{\epsilon}(X)$.

\begin{figure}[!ht]
\centering
\begin{tabular}{ll}
\subfloat[$\re U_1^{\epsilon}$]{
\includegraphics[height=.35\textwidth]{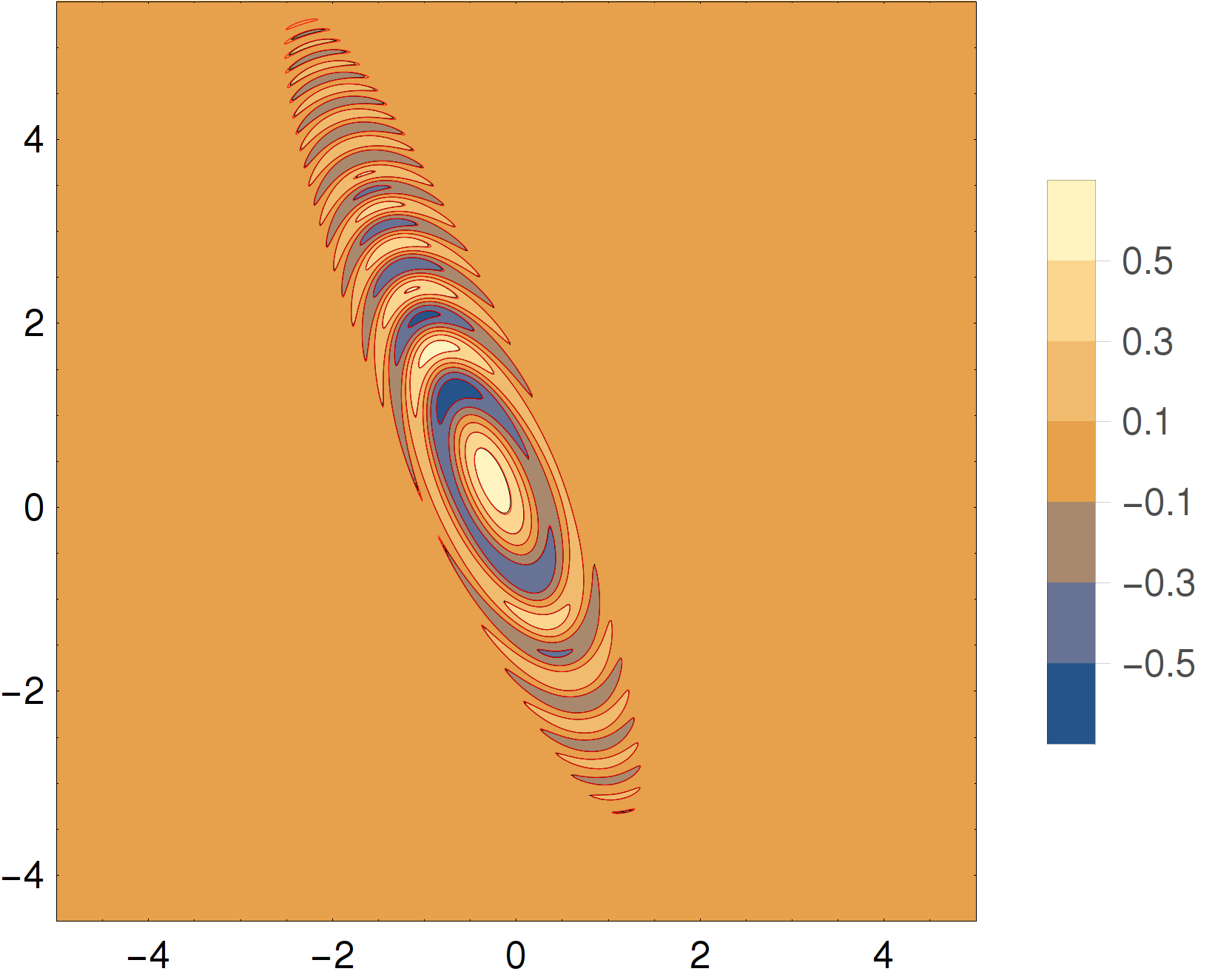}
} &
\subfloat[$\im U_1^{\epsilon}$]{
\includegraphics[height=.35\textwidth]{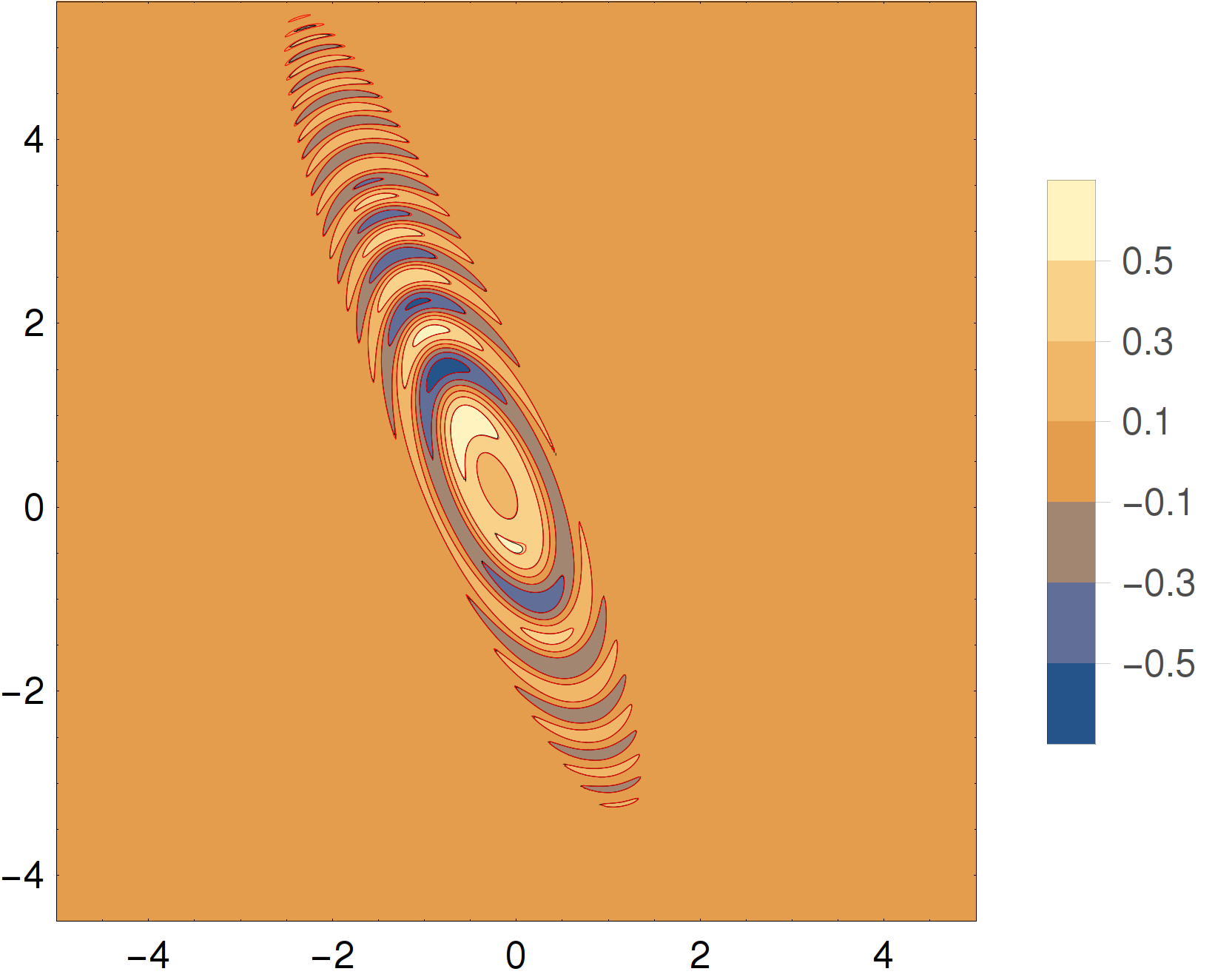}
} \\[30pt]
\subfloat[$\re U_2^{\epsilon}$]{
\includegraphics[height=.35\textwidth]{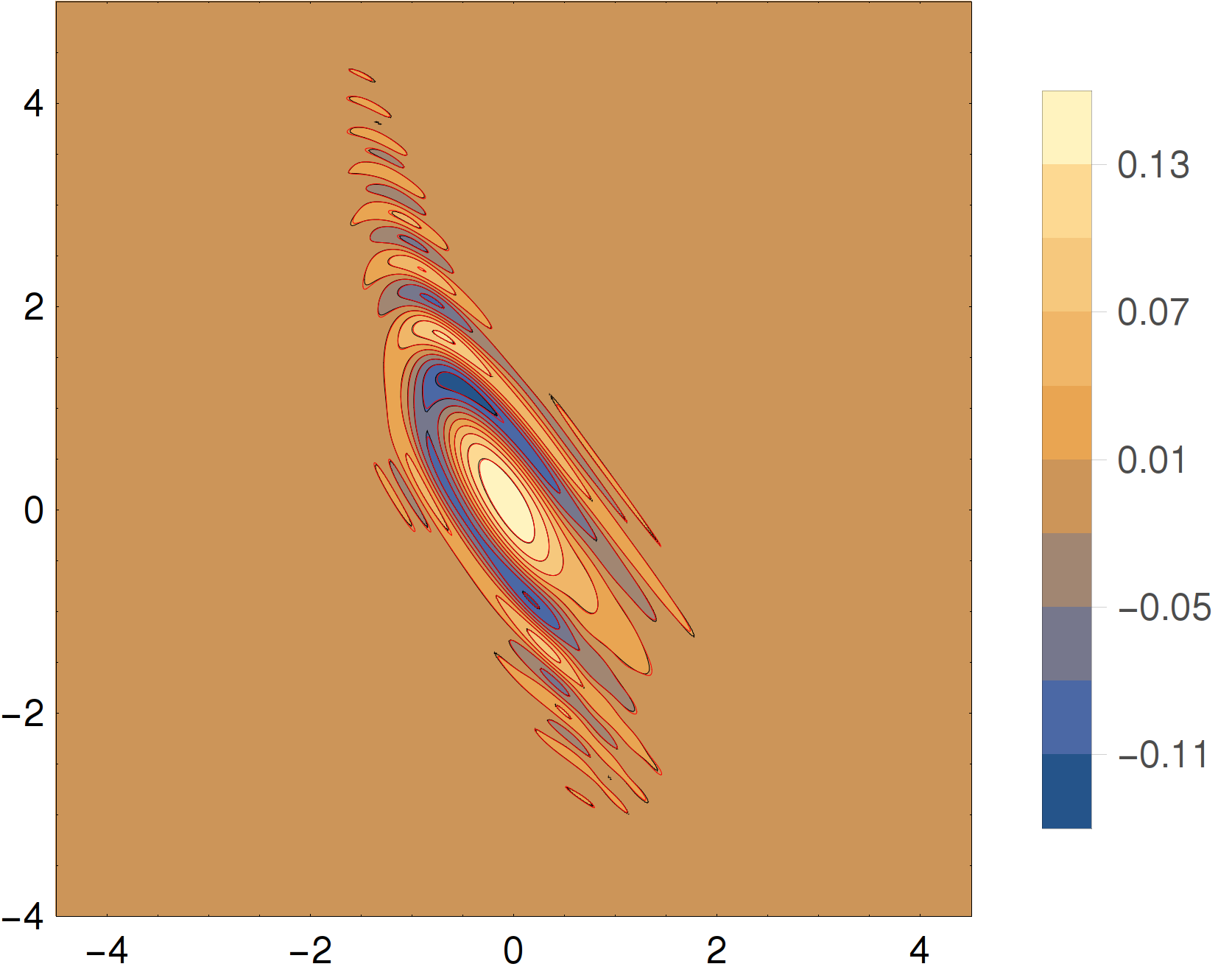}
} &
\subfloat[$\im U_2^{\epsilon}$]{
\includegraphics[height=.35\textwidth]{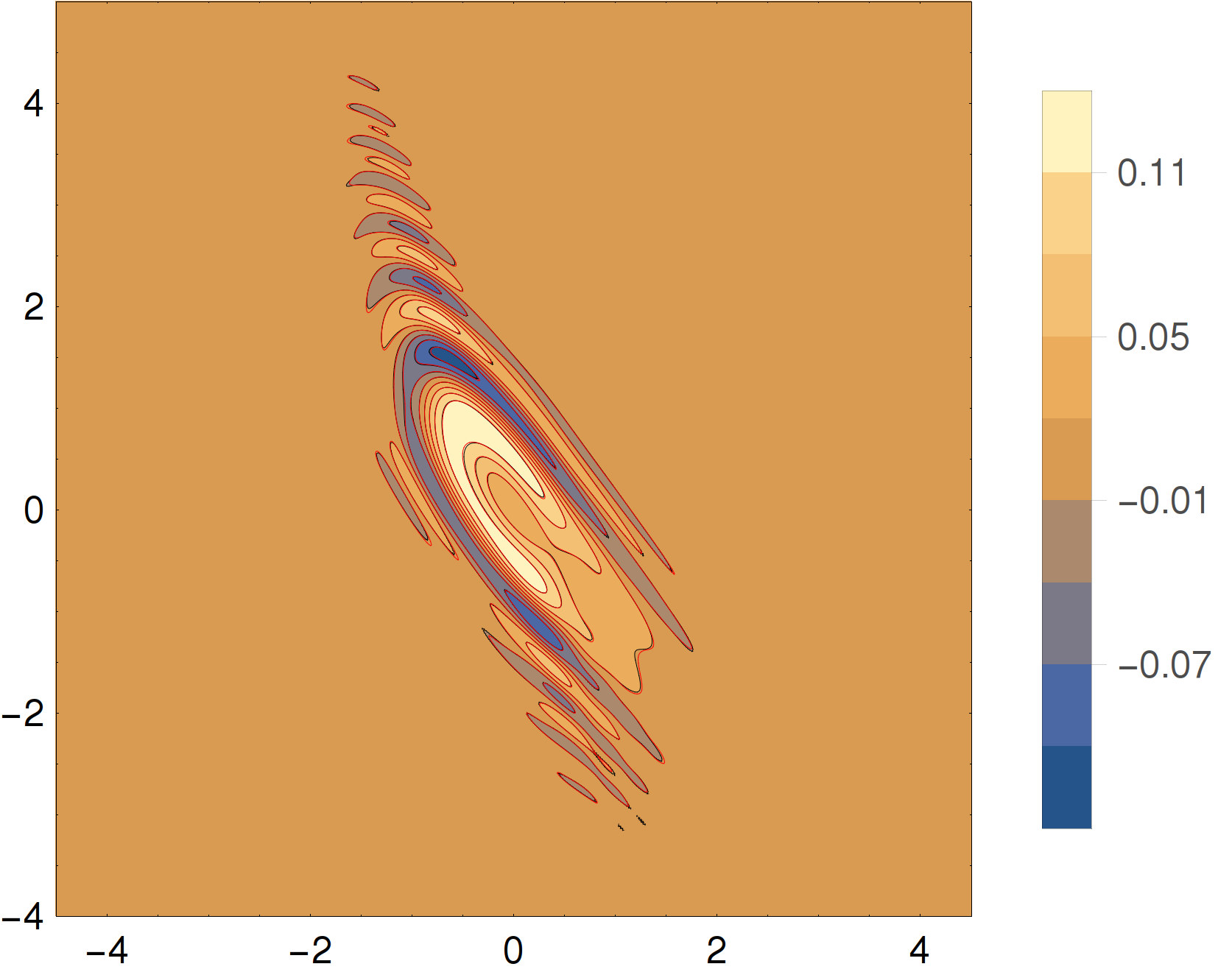}
}
\end{tabular}
\caption{Comparison between the numerical result and the reference
solution. The white lines show the contours of the reference solution,
and the color shading and black lines (almost covered by the red
lines) give the contours of the solution of the surface hopping
Gaussian beam method.}
\label{fig:ex1_comp}
\end{figure}

In this example, the coefficients $\alpha_i$ are linear in $x$,
$\beta_i$ are quadratic, and $\gamma_{ij}$ are constants. Therefore
the $O(\epsilon^{\frac{1}{2}})$ residuals in \eqref{eq:G_eq} and
\eqref{eq:I} vanish. Hence the expectation $U_i^{\epsilon}(t,x)$ of
the surface hopping Gaussian beam method is the exact solution of
\eqref{eq:ex1}. Since \eqref{eq:ex1} does not have an explicit
analytical solution, we solve the system with a finite volume method
and the result is used a reference solution. The scheme is third order
in both time and space, and the mesh is uniform with grid size $\Delta
x_1 = \Delta x_2 = 0.02$. Figure \ref{fig:ex1_comp} shows the
numerical results for $\epsilon = 0.1$ at $t = 0.5$, which are
obtained with 100,000 trajectories. The isolines are almost identical
to the reference solution, which indicates the feasibility of the
surface hopping Gaussian beam method.

\begin{figure}[!ht]
\centering
\subfloat[$\epsilon=0.1$]{
\includegraphics[height=.28\textwidth]{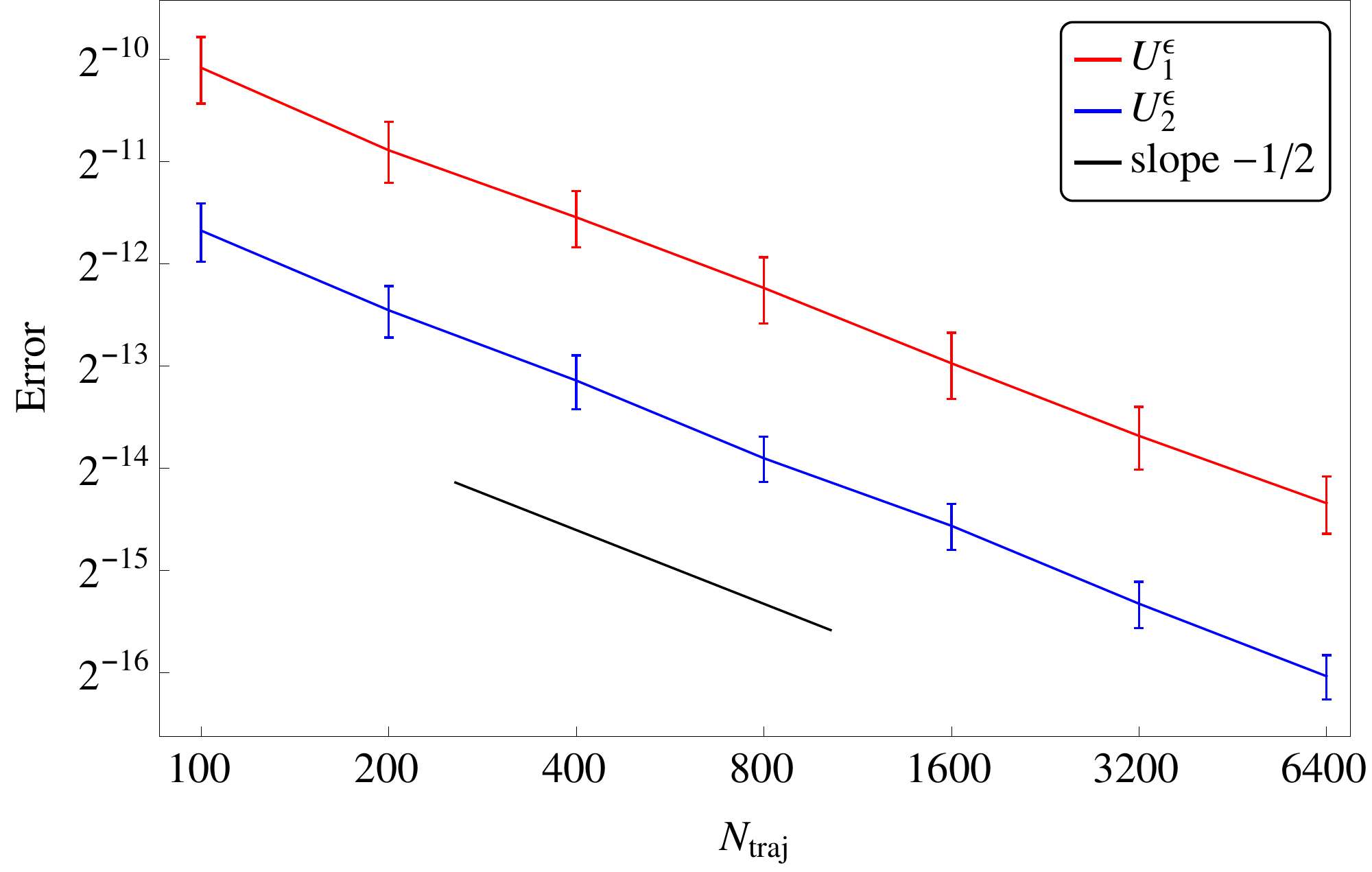}
} \hspace{20pt}
\subfloat[$\epsilon=0.5$]{
\includegraphics[height=.28\textwidth]{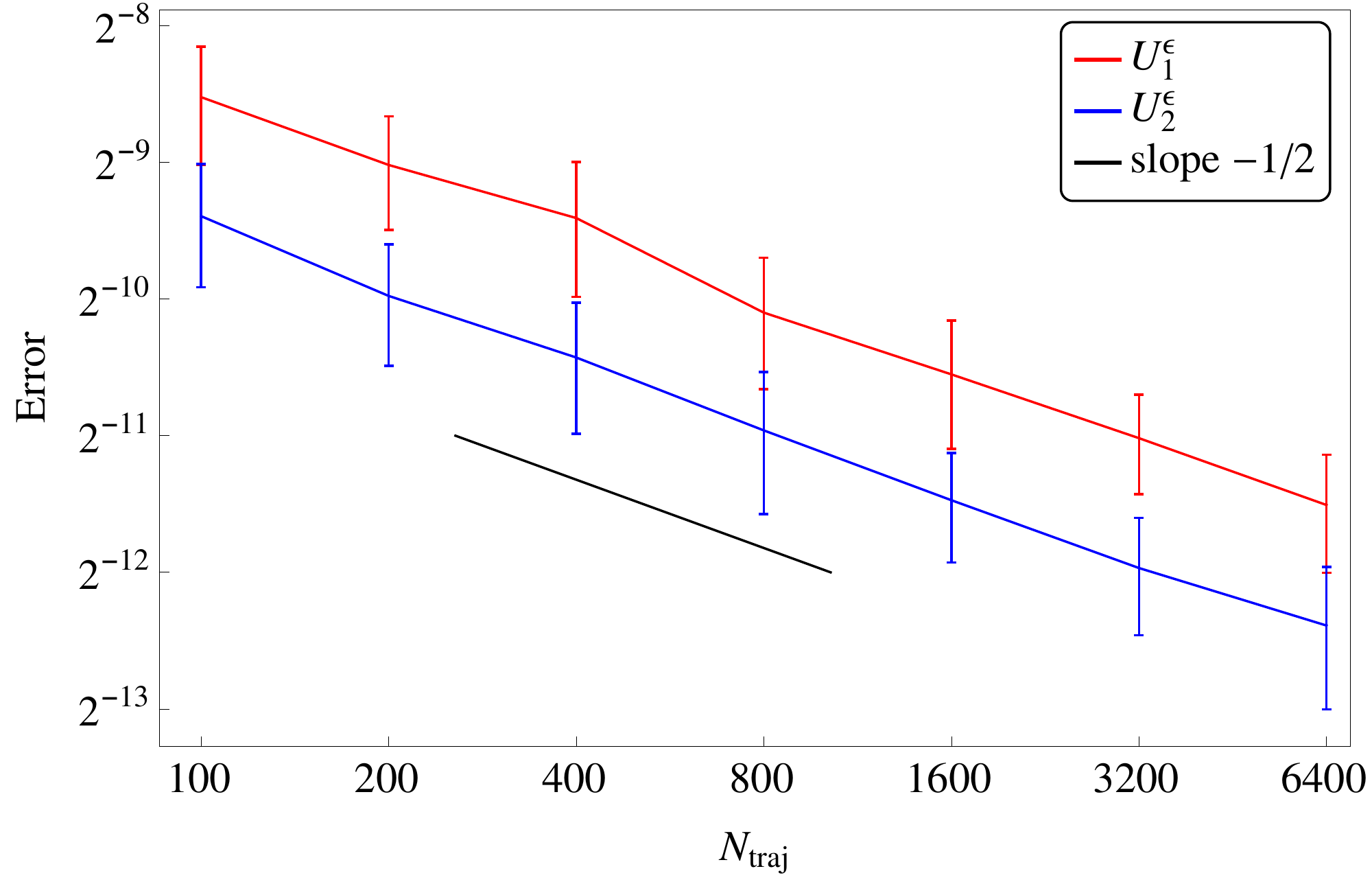}
}
\caption{Log-log plot of the $L^2$ errors for $U_1^{\epsilon}$ and
$U_2^{\epsilon}$ with different number of trajectories}
\label{fig:ex1_error}
\end{figure}

Below we study the convergence of the surface hopping Gaussian beam
method with respect to the number of trajectories $N_{\mathrm{traj}}$.
We run the above simulation with 100, 200, 400, 800, 1600, 3200, 6400
trajectories and compute the $L^2$ difference between the numerical
solution and the reference solution. The result is plotted in Figure
\ref{fig:ex1_error}, where the case $\epsilon = 0.5$ is also included.
The errors in the figure are obtained by averaging the erros of 100
runs in each case, and the error bars show the standard deviations.
The figure clearly shows that both the error and its standard
deviation decreases as $O(N_{\mathrm{traj}}^{-1/2})$ when the number
of trajectories increases, which agrees with the general property of
stochastic numerical methods.

\subsection{An example with nonlinear coefficients}
In this section, we consider the system
\begin{equation} \label{eq:ex2}
\begin{aligned}
& \frac{\partial u_1^{\epsilon}}{\partial t} -
  \sin(x_1 + x_2) \frac{\partial u_1^{\epsilon}}{\partial x_1} -
  \cos(x_1 - x_2) \frac{\partial u_1^{\epsilon}}{\partial x_2} =
\frac{\mathrm{i}}{\epsilon} \sin(|x|^2) u_1^{\epsilon} -
  |x|^2 (u_1^{\epsilon} - 5 u_2^{\epsilon}), \\
& \frac{\partial u_2^{\epsilon}}{\partial t} -
  \sin(x_1 - x_2) \frac{\partial u_2^{\epsilon}}{\partial x_1} -
  \cos(x_1 + x_2) \frac{\partial u_2^{\epsilon}}{\partial x_2} =
\frac{\mathrm{i}}{\epsilon} \cos(|x|^2) u_2^{\epsilon} -
  |x|^2 \left( 5 u_1^{\epsilon} + \frac{1}{2} u_2^{\epsilon} \right),
\end{aligned}
\end{equation}
with initial conditions
\begin{equation}
u_1^{\epsilon}(0,x) = \exp \left( -\frac{|x|^2}{2\epsilon} \right),
\qquad u_2^{\epsilon}(0,x) = 0.
\end{equation}
In this example, the coefficients are nonlinear in $x$, and therefore
the $O(\epsilon^{\frac{1}{2}})$ residuals in \eqref{eq:G_eq} and
\eqref{eq:I} exist, which will lead to some discrepancy between
$u_i^{\epsilon}$ and $U_i^{\epsilon}$. Here we again use the finite
volume solver to provide a reference solution, and the grid size is
$\Delta x_1 = \Delta x_2 = 0.005$.

\begin{figure}[!ht]
\centering
\subfloat[$U_1^{\epsilon}$ and the reference solution]{
\includegraphics[height=.33\textwidth]{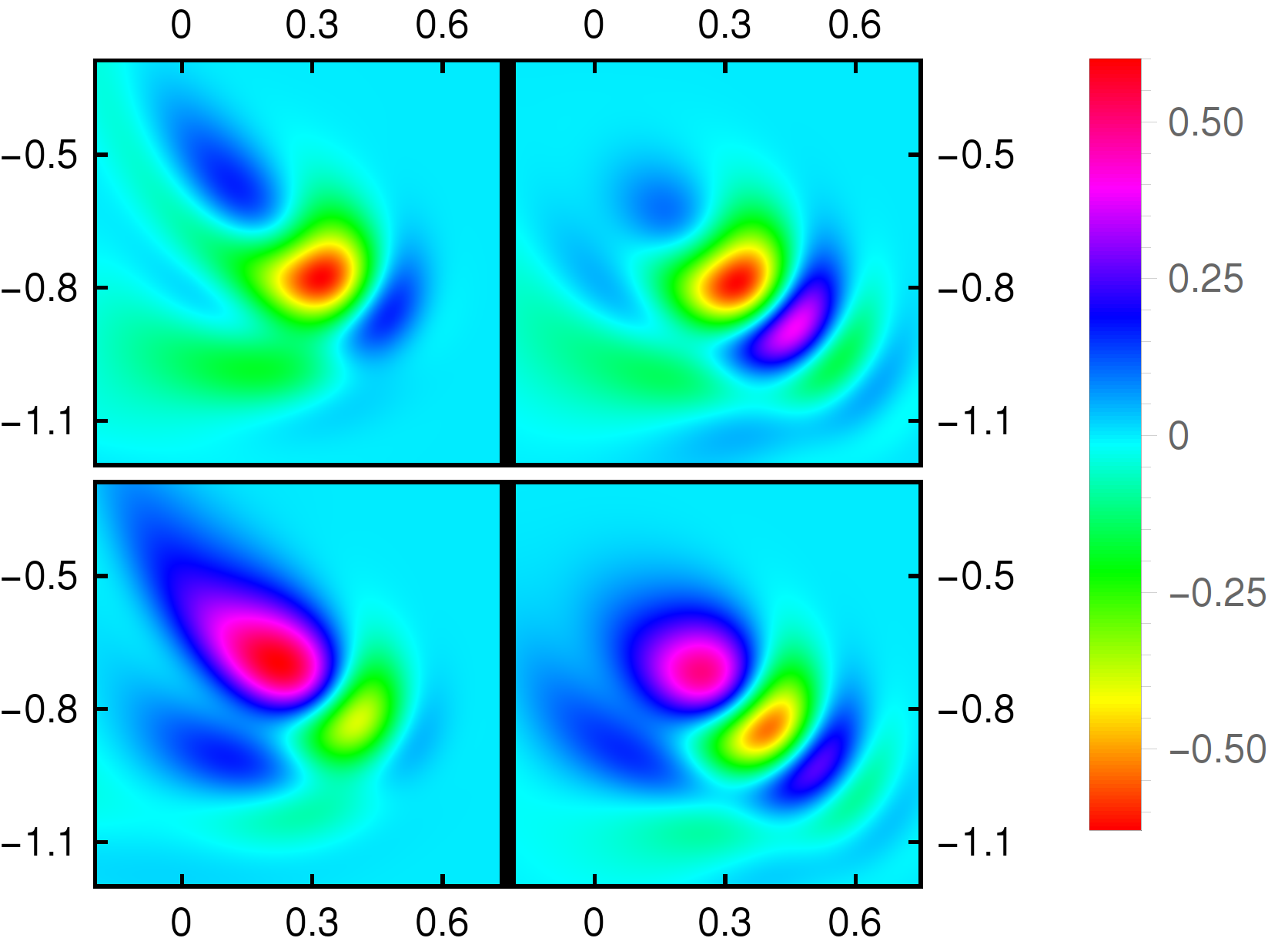}
} \hspace{20pt}
\subfloat[$U_2^{\epsilon}$ and the reference solution]{
\includegraphics[height=.33\textwidth]{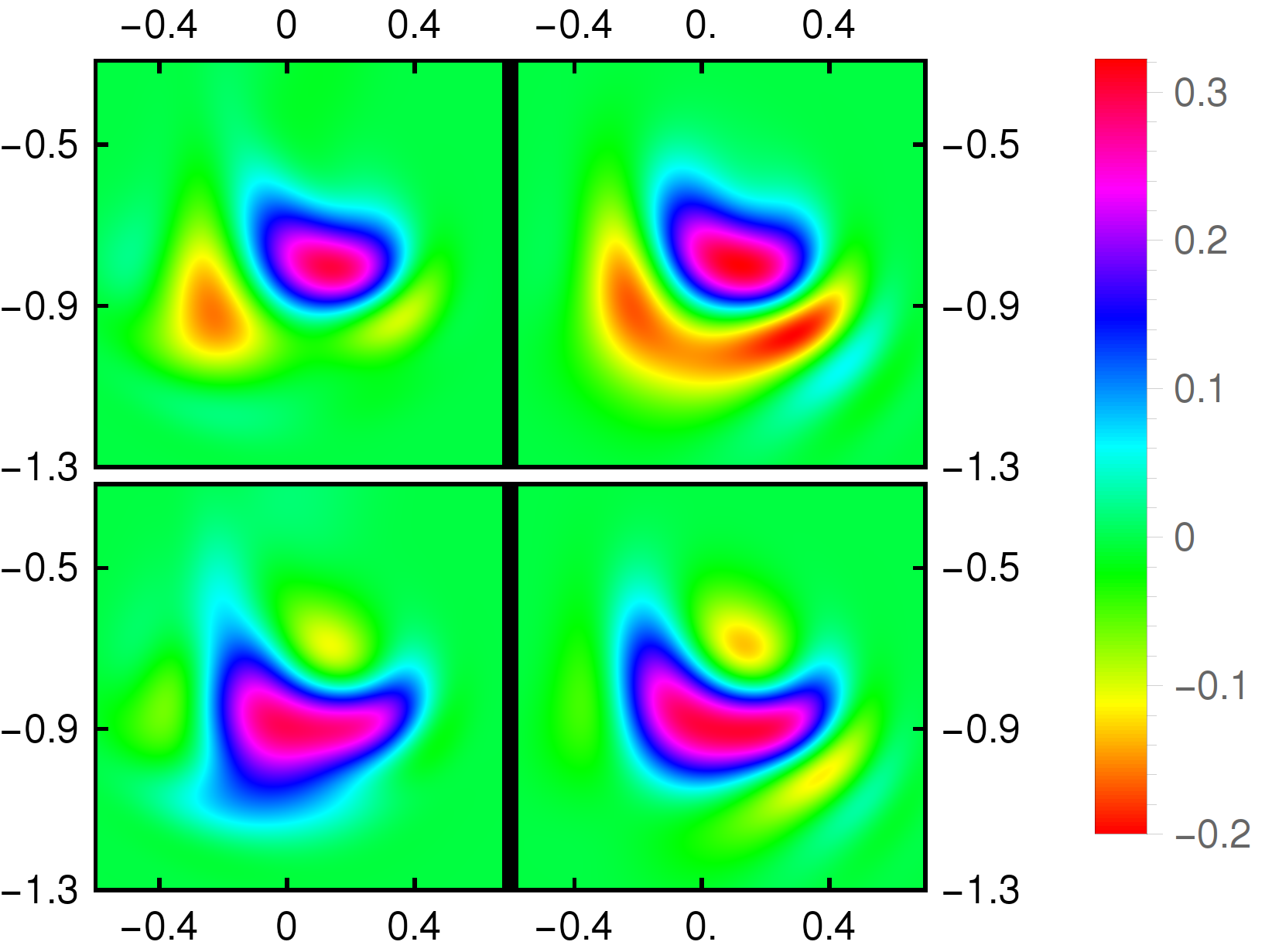}
}
\caption{Comparison between the numerical results and the reference
solutions with $\epsilon = 0.08$. The upper half of the figure shows
the real part of the solution, and the lower half shows the imaginary
part. The results of the surface hopping Gaussian beam method are
given to the right of the bold black bars, and the left parts are the
reference solutions.}
\label{fig:ex2_eps0p08_comp}
\end{figure}

\begin{figure}[!ht]
\centering
\subfloat[$U_1^{\epsilon}$ and the reference solution]{
\includegraphics[height=.33\textwidth]{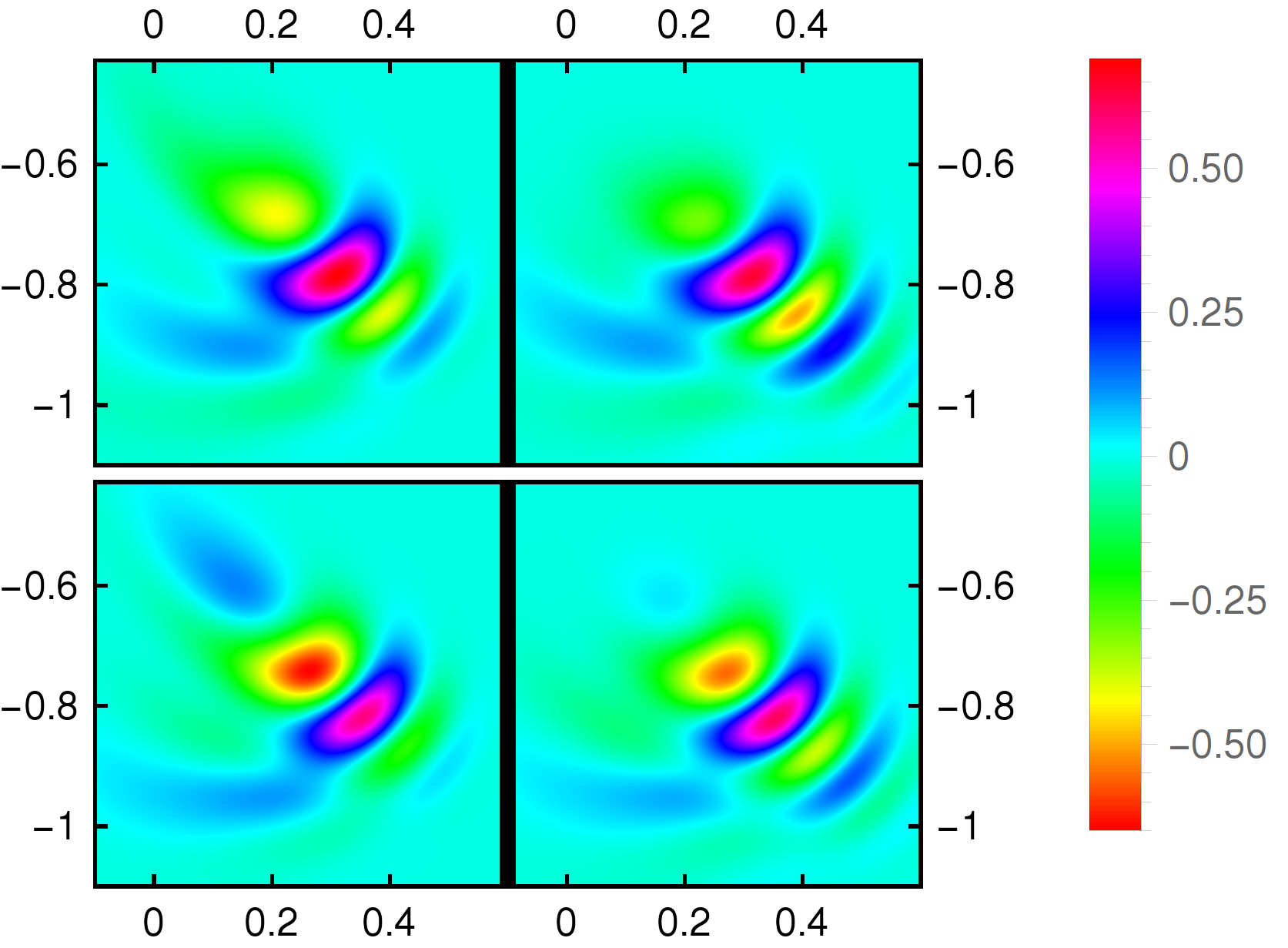}
} \hspace{20pt}
\subfloat[$U_2^{\epsilon}$ and the reference solution]{
\includegraphics[height=.33\textwidth]{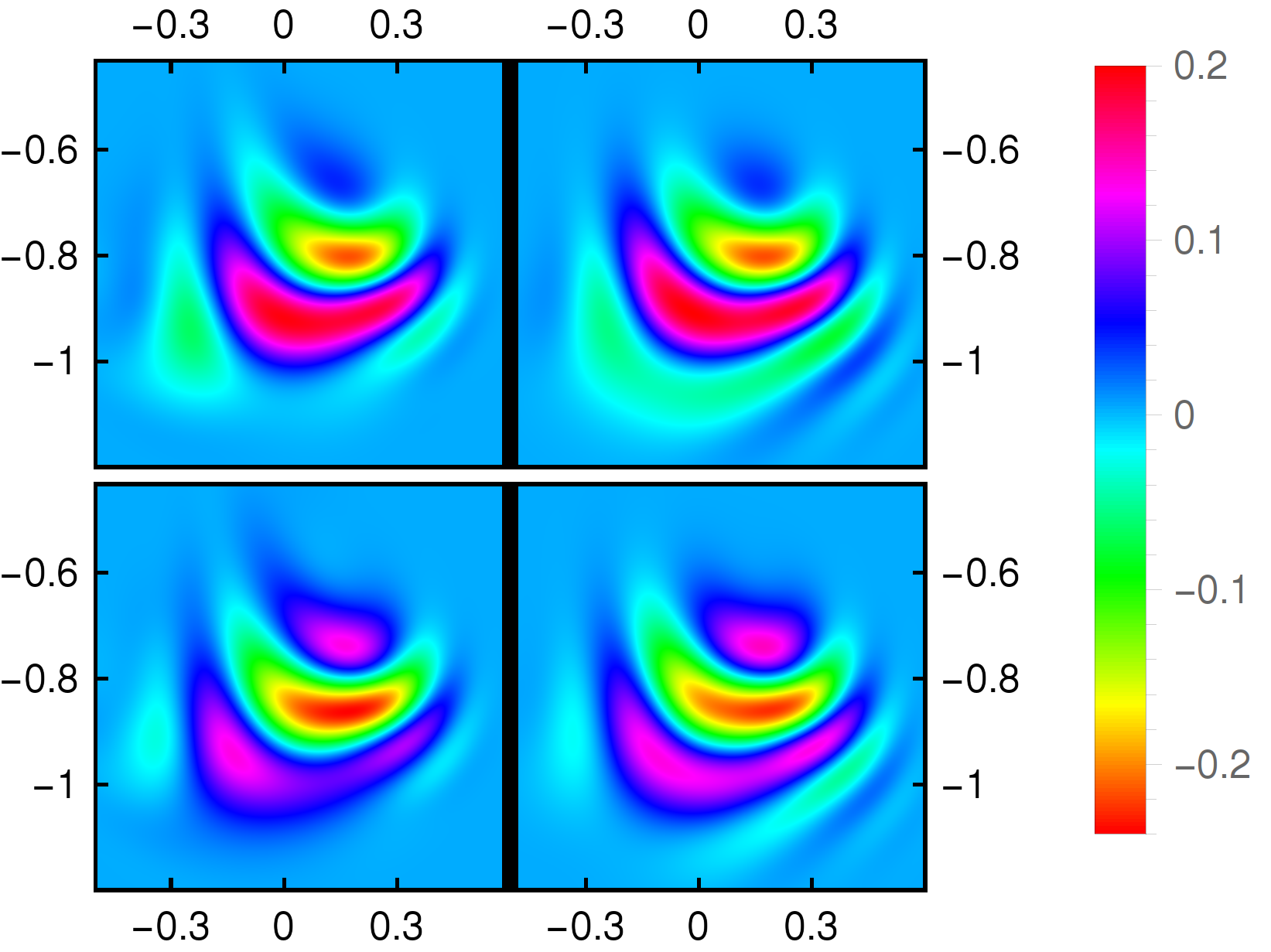}
}
\caption{Comparison between the numerical results and the reference
solutions with $\epsilon = 0.04$. See Figure
\ref{fig:ex2_eps0p08_comp} for the details.}
\label{fig:ex2_eps0p04_comp}
\end{figure}

\begin{figure}[!ht]
\centering
\subfloat[$U_1^{\epsilon}$ and the reference solution]{
\includegraphics[height=.33\textwidth]{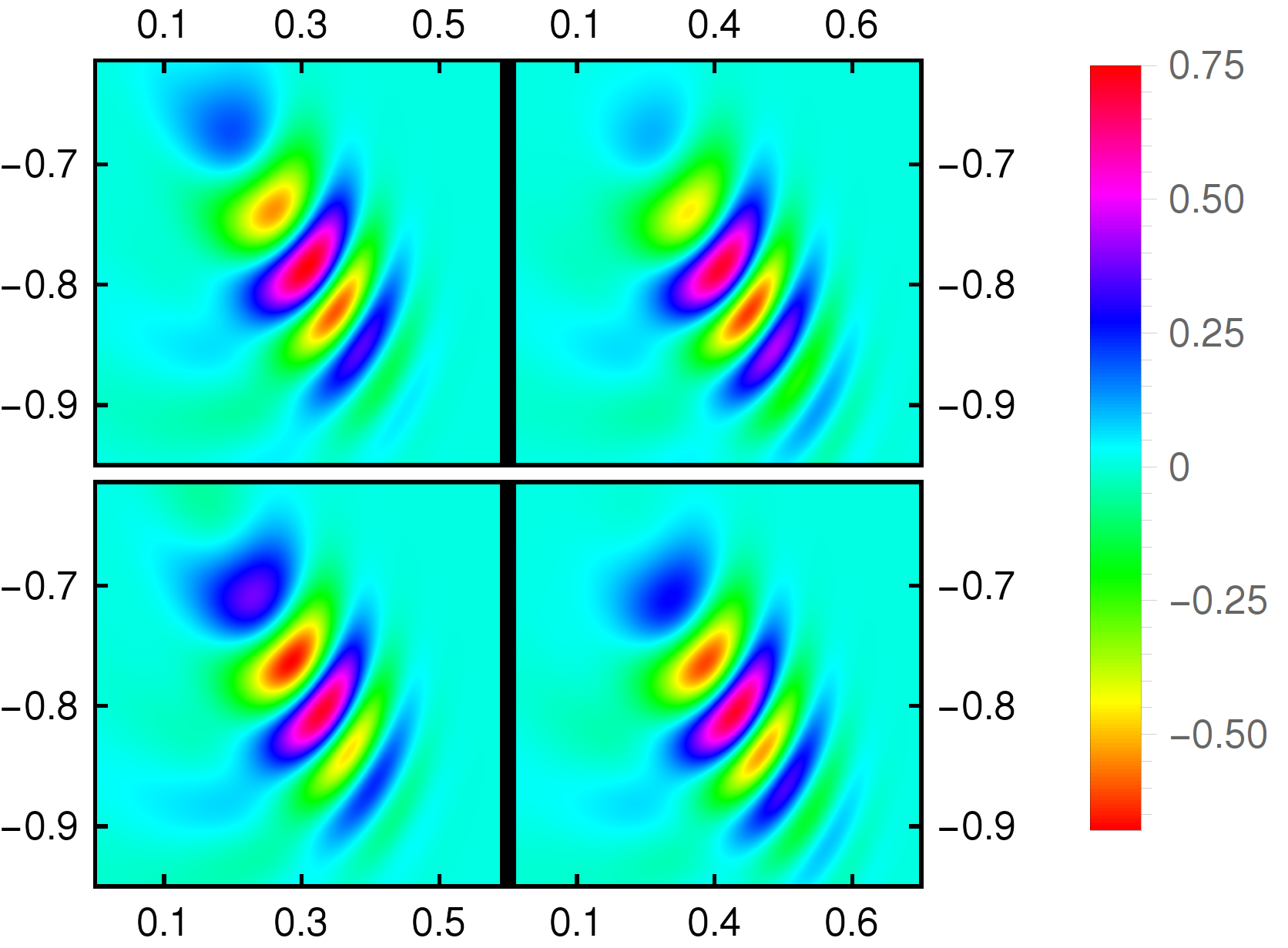}
} \hspace{20pt}
\subfloat[$U_2^{\epsilon}$ and the reference solution]{
\includegraphics[height=.33\textwidth]{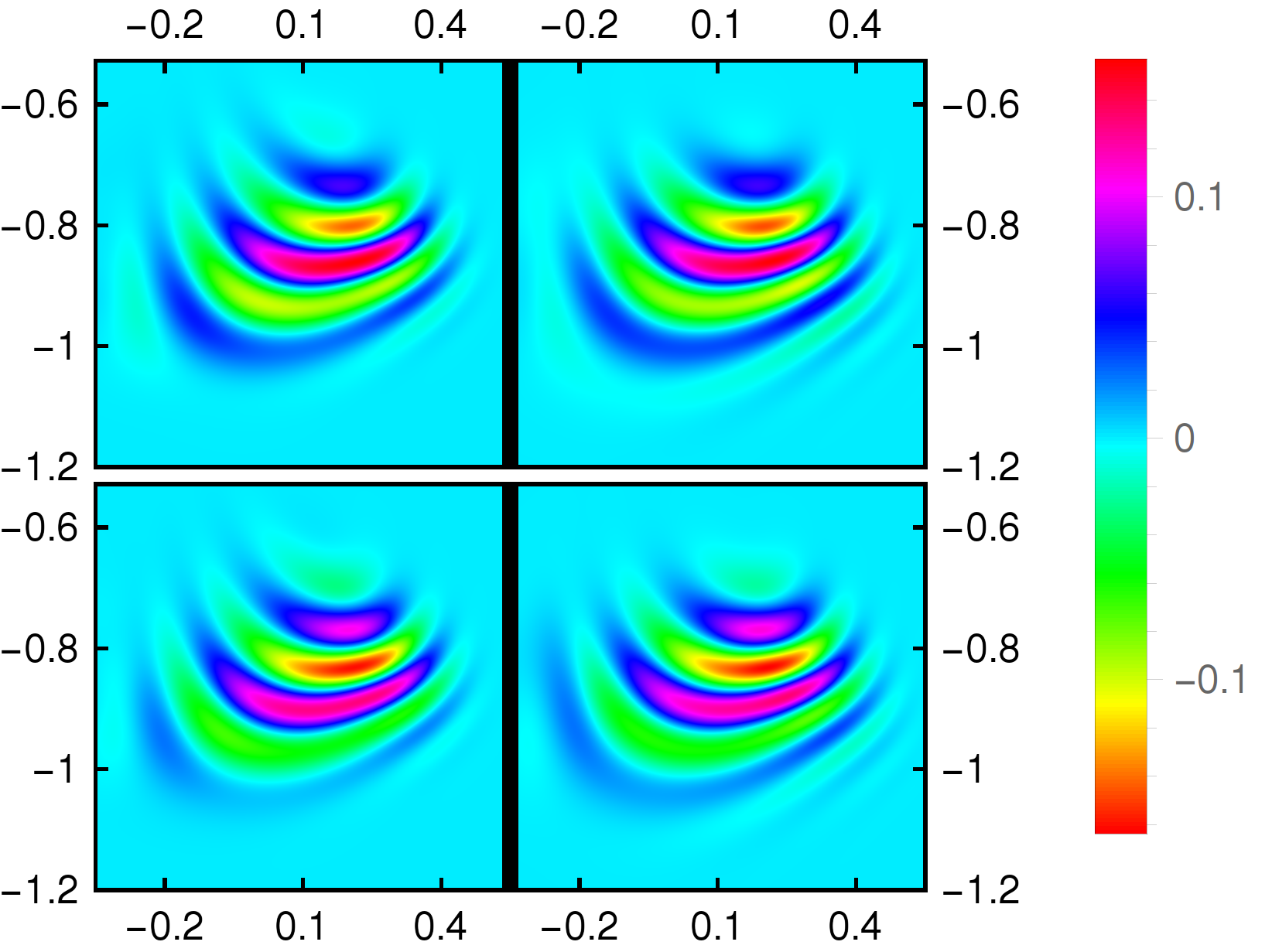}
}
\caption{Comparison between the numerical results and the reference
solutions with $\epsilon = 0.02$. See Figure
\ref{fig:ex2_eps0p08_comp} for the details.}
\label{fig:ex2_eps0p02_comp}
\end{figure}

\begin{figure}[!ht]
\centering
\subfloat[$U_1^{\epsilon}$ and the reference solution]{
\includegraphics[height=.33\textwidth]{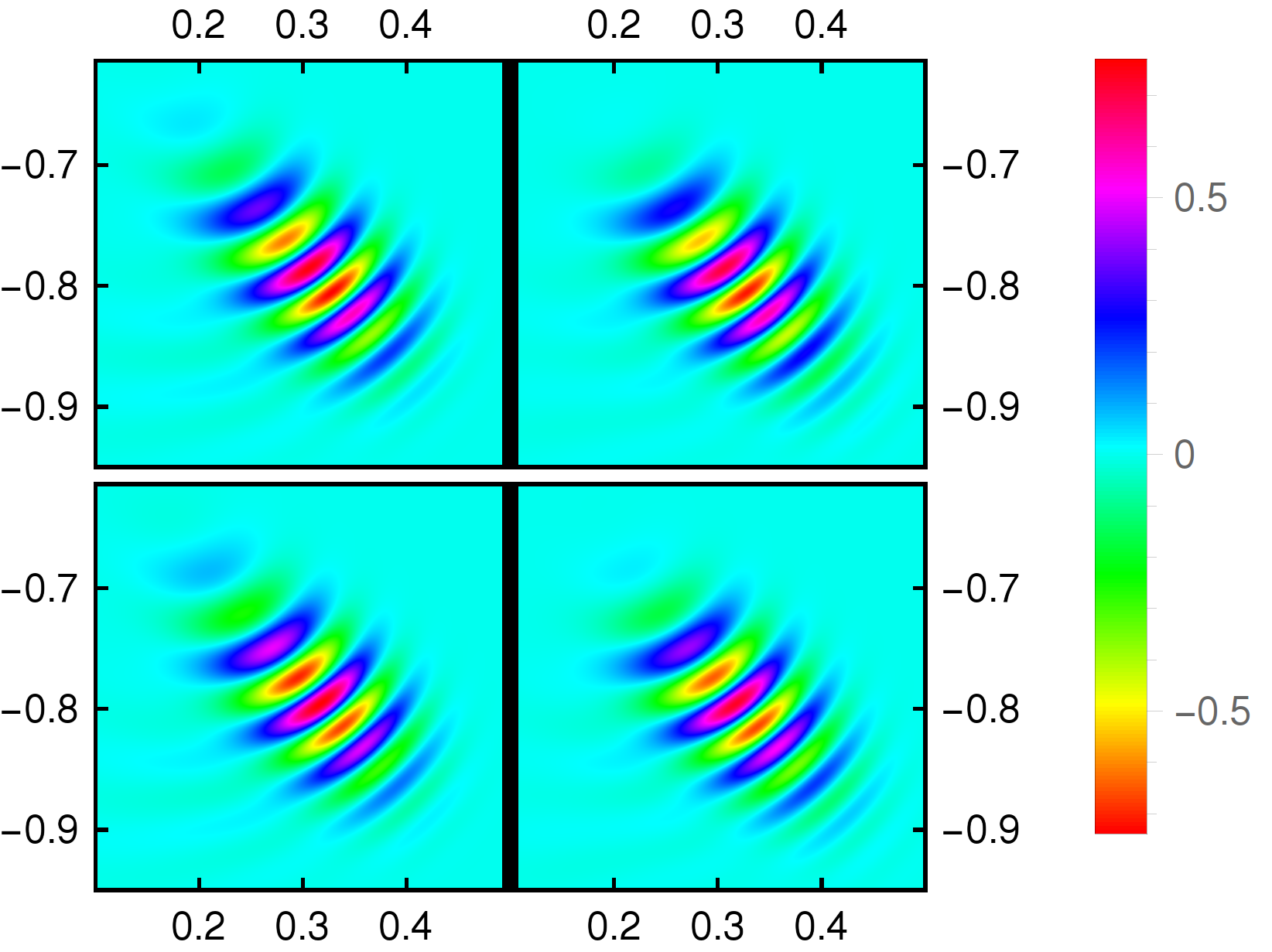}
} \hspace{20pt}
\subfloat[$U_2^{\epsilon}$ and the reference solution]{
\includegraphics[height=.33\textwidth]{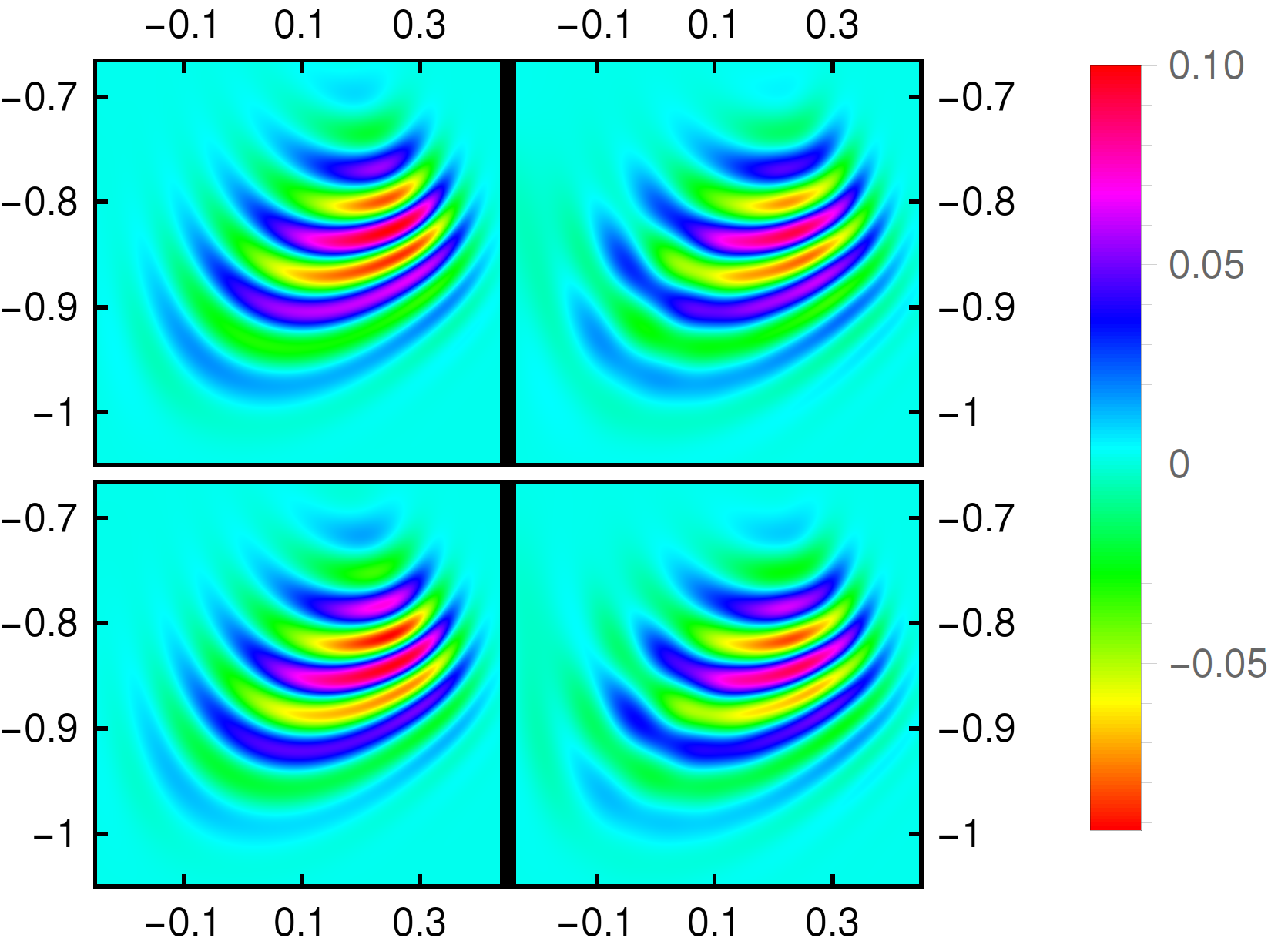}
}
\caption{Comparison between the numerical results and the reference
solutions with $\epsilon = 0.01$. See Figure
\ref{fig:ex2_eps0p08_comp} for the details.}
\label{fig:ex2_eps0p01_comp}
\end{figure}

We consider four different values of $\epsilon$. Figure
\ref{fig:ex2_eps0p08_comp}--\ref{fig:ex2_eps0p01_comp}
respectively show the results for $\epsilon = 0.08, 0.04, 0.02, 0.01$.
All the simulations are performed using 10,000 trajectories. As
expected, the solution spreads wider when $\epsilon$ increases, and
the oscillation gets stronger when $\epsilon$ decreases. In all the
cases, our method generates reasonable approximate solutoins, while it
can still be observed from the figures that better approximation of
the reference solution is obtained for smaller $\epsilon$.

\subsection{Quantum-classical Liouville equation}
In this section, we consider the application of the surface hopping
Gaussian beam method to the quantum-classical Liouville equation
\cite{Schutte1999}, which describes the composed system with heavy and
light particles. The quantum-classical Liouville equation (QCLE) has
exactly the form \eqref{eq:extended}, and thus our algorithm can be
directly applied. Here we study the two level system as in
\cite{Chai2015}, and write the equations as
\begin{equation} \label{eq:QCLE}
\begin{aligned}
\frac{\partial u_{11}^{\epsilon}(t,r,p)}{\partial t} &=
  -p \cdot \nabla_r u_{11}^{\epsilon}(t,r,p) +
  \nabla_r E_1(r) \cdot \nabla_p u_{11}^{\epsilon}(t,r,p) +
  [p \cdot \overline{d_{21}(r)}] u_{21}^{\epsilon}(t,r,p) +
  [p \cdot d_{21}(r)] \overline{u_{21}^{\epsilon}(t,r,p)}, \\
\frac{\partial u_{22}^{\epsilon}(t,r,p)}{\partial t} &=
  -p \cdot \nabla_r u_{22}^{\epsilon}(t,r,p) +
  \nabla_r E_2(r) \cdot \nabla_p u_{11}^{\epsilon}(t,r,p) -
  [p \cdot \overline{d_{21}(r)}] u_{21}^{\epsilon}(t,r,p) -
  [p \cdot d_{21}(r)] \overline{u_{21}^{\epsilon}(t,r,p)}, \\
\frac{\partial u_{21}^{\epsilon}(t,r,p)}{\partial t} &=
  -\frac{\mathrm{i}}{\epsilon} [E_2(r)-E_1(r)] u_{21}^{\epsilon}(t,r,p)
  -p \cdot \nabla_r u_{21}^{\epsilon}(t,r,p) +
  \frac{1}{2} \nabla_r [E_1(r) + E_2(r)]
    \cdot \nabla_p u_{21}^{\epsilon}(t,r,p) \\
& \qquad + [p \cdot (d_{11}(r) - d_{22}(r))] u_{21}^{\epsilon}(t,r,p)
  - [p \cdot d_{21}(r)] u_{11}^{\epsilon}(t,r,p)
  + [p \cdot d_{21}(r)] u_{22}^{\epsilon}(t,r,p),
\end{aligned}
\end{equation}
where $E_1(r), E_2(r) \in \mathbb{R}$ are adiabatic energy surfaces,
$r \in \mathbb{R}^n$ and $p \in \mathbb{R}^n$ stand respectivly for the
position and momentum, and $d_{ij}(r) \in \mathbb{R}^n$ is used to
denote the first-order adiabatic operator. We refer the readers to
\cite{Schutte1999,Chai2015} for the derivation of the above equations,
and here we just mention that $E_i(r)$ and $d_{ij}(r)$ can be derived
from the potential energy matrix acting on the heavy particles $V(r)
\in \mathbb{C}^{2\times 2}$. Assuming $V(r)$ to be real and symmetric,
we have
\begin{equation}
\begin{gathered}
E_{1,2}(r) = \frac{1}{2} [V_{11}(r) + V_{22}(r)] \pm
  \frac{1}{2} \sqrt{[V_{11}(r) - V_{22}(r)]^2 + 4V_{12}^2(r)}, \\
d_{21}(r) = \frac{[V_{11}(r) - V_{22}(r)] \nabla_r V_{12}(r) -
    V_{12}(r) \nabla_r[V_{11}(r) - V_{22}(r)]}
  {[V_{11}(r) - V_{22}(r)]^2 + 4 V_{12}^2(r)}, \quad
d_{11}(r) = d_{22}(r) = 0.
\end{gathered}
\end{equation}
The equations \eqref{eq:QCLE} are slightly different from the more
conventional QCLE derived in \cite{Kapral1999,Horenko2002}. The system
\eqref{eq:QCLE} can be considered as the simplified version of the
quantum-classical Liouville equations, where a momentum-jump term is
absent. The simplification yields a $O(\epsilon)$ approximation to the
quantum-classical expectation values of physical quantities
\cite{Schutte1999}. We refer the readers to \cite{Horenko2002,
Ryabinkin2014, Cai2017} for the comparisons of the two models.

Below we will present two examples to test our method.

\subsubsection{Extended coupling with reflection}
\begin{figure}[!ht]
\centering
\subfloat[$\epsilon = 1/16$]{
\includegraphics[width=.28\textwidth]{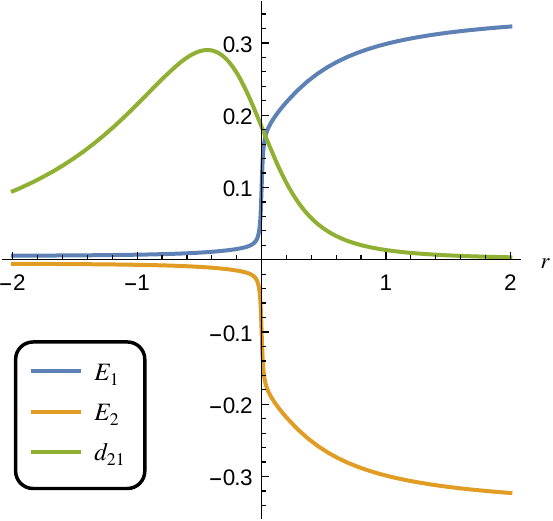}
} \hspace{20pt}
\subfloat[$\epsilon = 1/32$]{
\includegraphics[width=.28\textwidth]{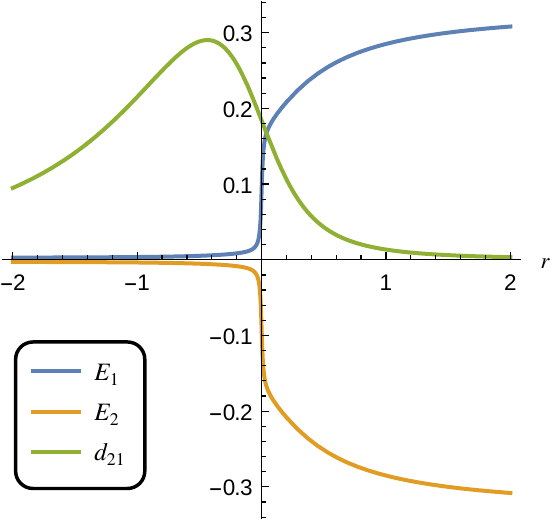}
} \hspace{20pt}
\subfloat[$\epsilon = 1/64$]{
\includegraphics[width=.28\textwidth]{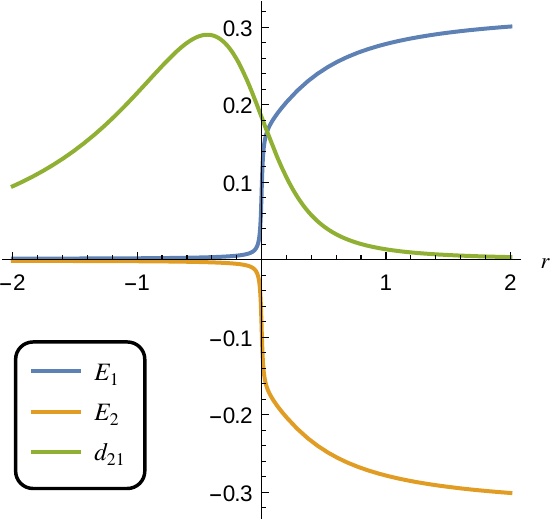}
}
\caption{Funtions $E_1(r)$, $E_2(r)$ and $d_{21}(r)$ for various
$\epsilon$}
\label{fig:ECR_E}
\end{figure}

This example is from one of Tully's examples in \cite{Tully1990}. The
number of dimensions $n = 1$, and the matrix $V(r)$ is given as
\begin{equation}
V(r) = \begin{pmatrix}
  V_{11}(r) & V_{12}(r) \\ V_{12}(r) & V_{22}(r)
\end{pmatrix} = F^{\delta}(r) \begin{pmatrix}
  \frac{1}{20} &
    \frac{1}{10} \left( \arctan(2r) + \frac{\pi}{2} \right) \\
  \frac{1}{10} \left( \arctan(2r) + \frac{\pi}{2} \right) &
    -\frac{1}{20}
\end{pmatrix},
\end{equation}
where $F^{\delta}(r) = \frac{1}{\pi} \left( \arctan(100r) +
\frac{\pi}{2} + \delta \right)$. Here we choose $\delta = 5 \epsilon$
and the initial condition
\begin{equation}
u_{11}^{\epsilon}(0,x) = u_{12}^{\epsilon}(0,x) = 0, \quad
u_{22}^{\epsilon}(0,x) = \frac{1}{\sqrt{32\pi} \epsilon} \exp \left(
  -\frac{(r-r_0)^2}{\epsilon} - \frac{(p-p_0)^2}{\epsilon}
\right),
\end{equation}
and the initial position and momentum are set to be $r_0 = -1.5$ and
$p_0 = 1.5$. Three cases $\epsilon = 1/16, 1/32, 1/64$ are considered,
and in these cases, the energy surfaces $E_{1,2}(r)$ and the function
$d_{21}(r)$ are plotted in Figure \ref{fig:ECR_E}. From
\eqref{eq:QCLE}, one can find that $E_{1,2}(r)$ describes the
``shape'' of the surfaces, and $d_{21}(r)$ determines the probability
of the surface hopping. In this example, $d_{21}(r)$ is indepedent of
$\epsilon$. If $r > 1.5$, there are hardly any hops between surfaces.

\begin{figure}[!ht]
\centering
\subfloat[$U_{11}^{\epsilon}$ and the reference solution]{
\includegraphics[height=.35\textwidth]{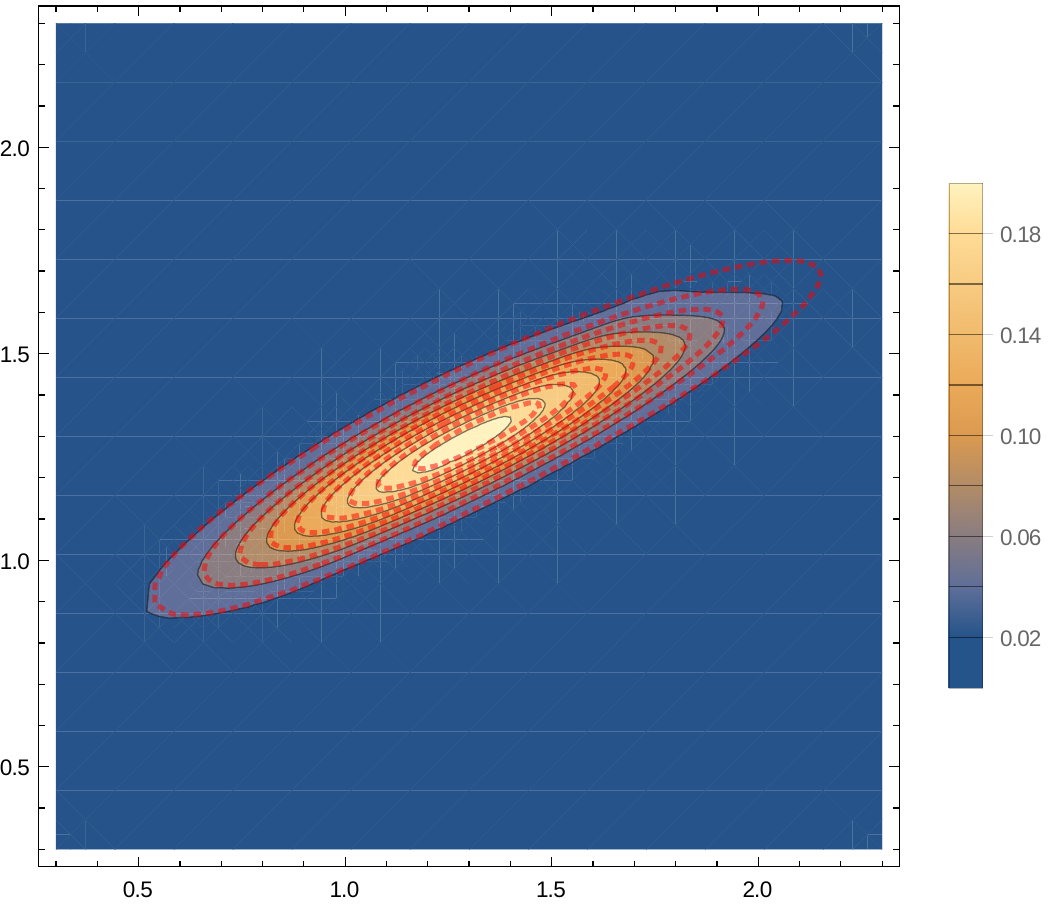}
} \hspace{20pt}
\subfloat[$U_{22}^{\epsilon}$ and the reference solution]{
\includegraphics[height=.35\textwidth]{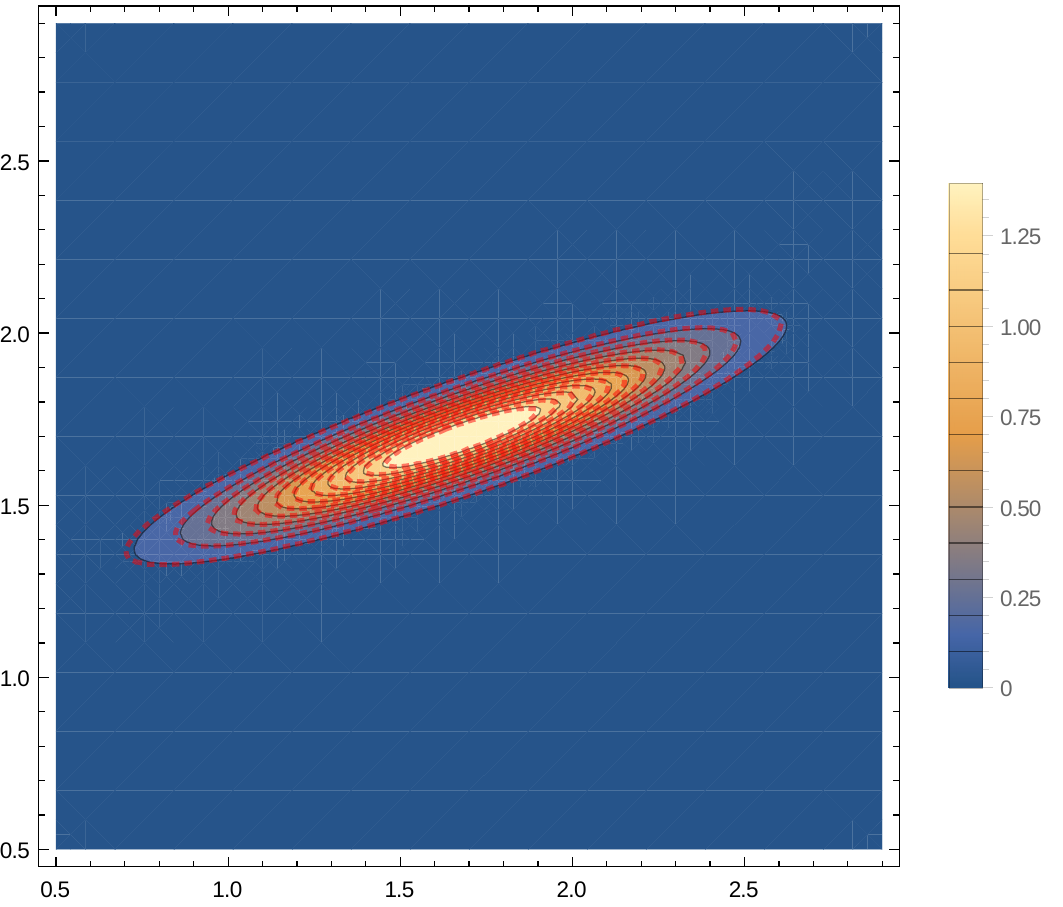}
}
\caption{Comparison between the numerical results (color shading and
black contour lines) and the reference solutions (red dotted contour
lines) for $\epsilon = 1/16$.}
\label{fig:eps_1d16}
\end{figure}

\begin{figure}[!ht]
\centering
\subfloat[$U_{11}^{\epsilon}$ and the reference solution]{
\includegraphics[height=.35\textwidth]{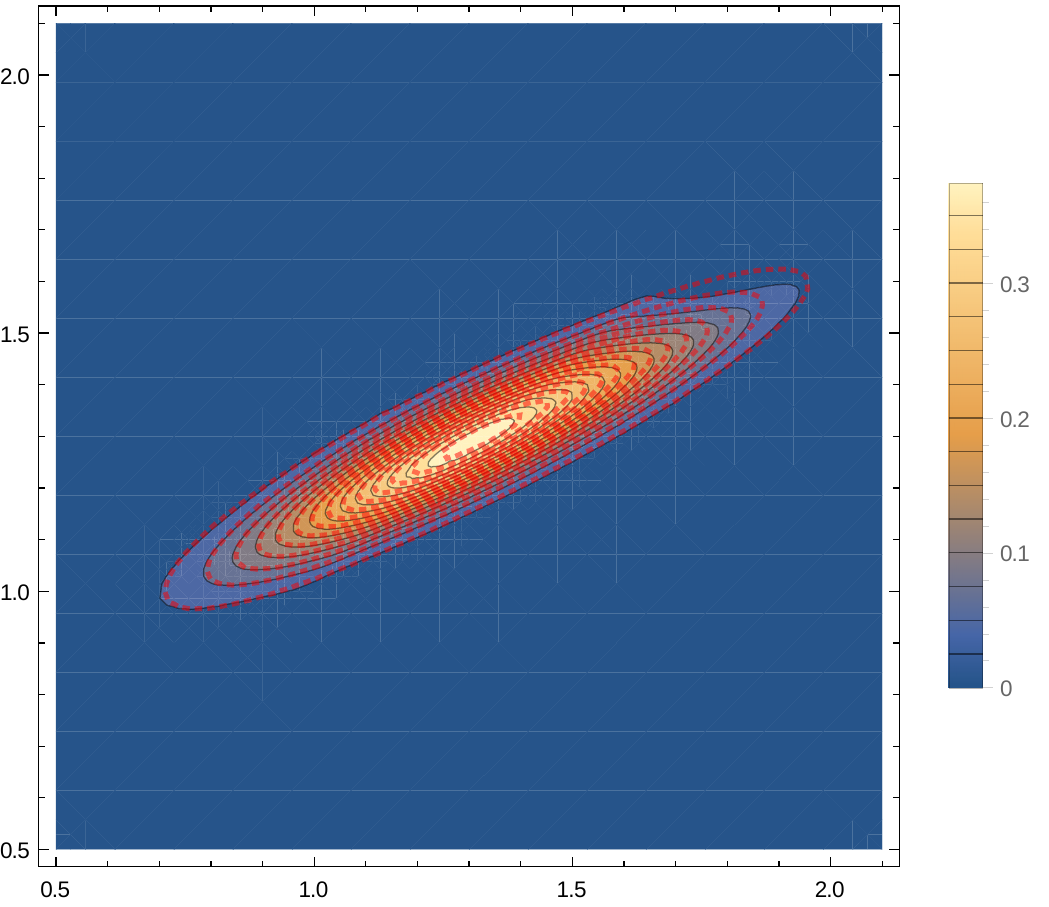}
} \hspace{20pt}
\subfloat[$U_{22}^{\epsilon}$ and the reference solution]{
\includegraphics[height=.35\textwidth]{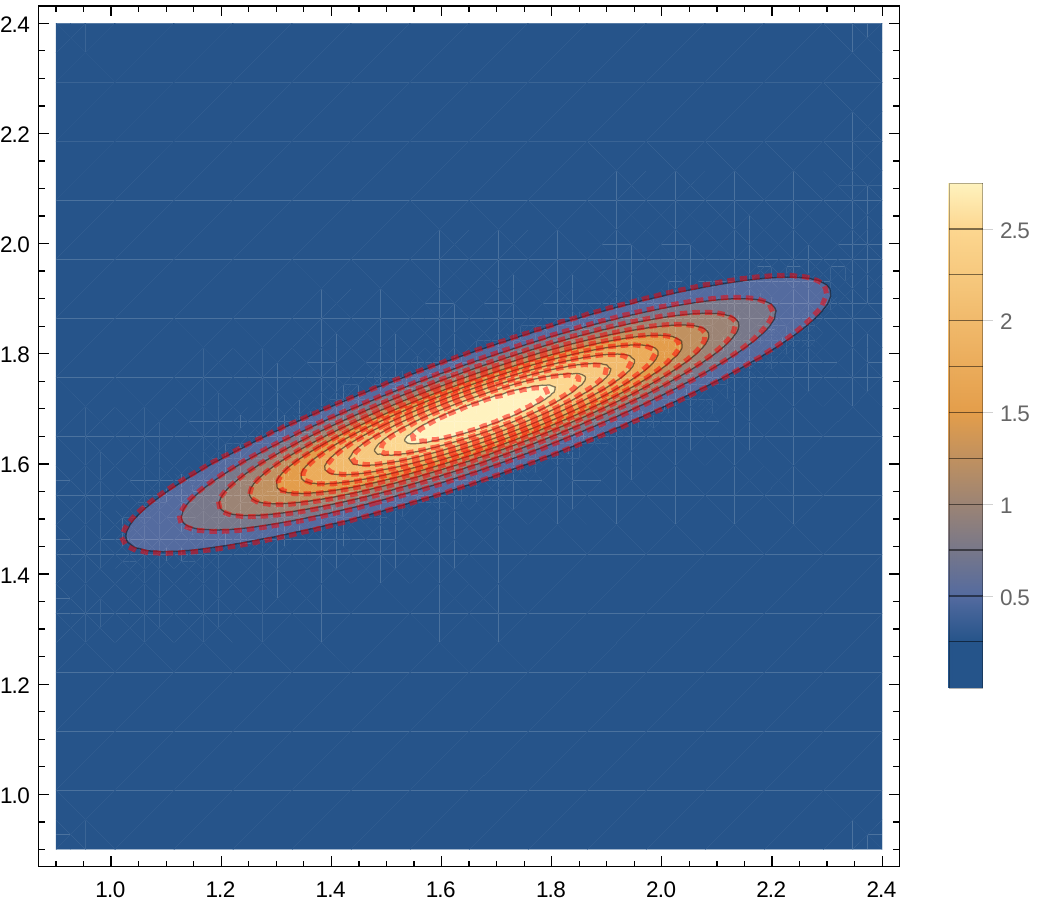}
}
\caption{Comparison between the numerical results (color shading and
black contour lines) and the reference solutions (red dotted contour
lines) for $\epsilon = 1/32$.}
\label{fig:eps_1d32}
\end{figure}

\begin{figure}[!ht]
\centering
\subfloat[$U_{11}^{\epsilon}$ and the reference solution]{
\includegraphics[height=.35\textwidth]{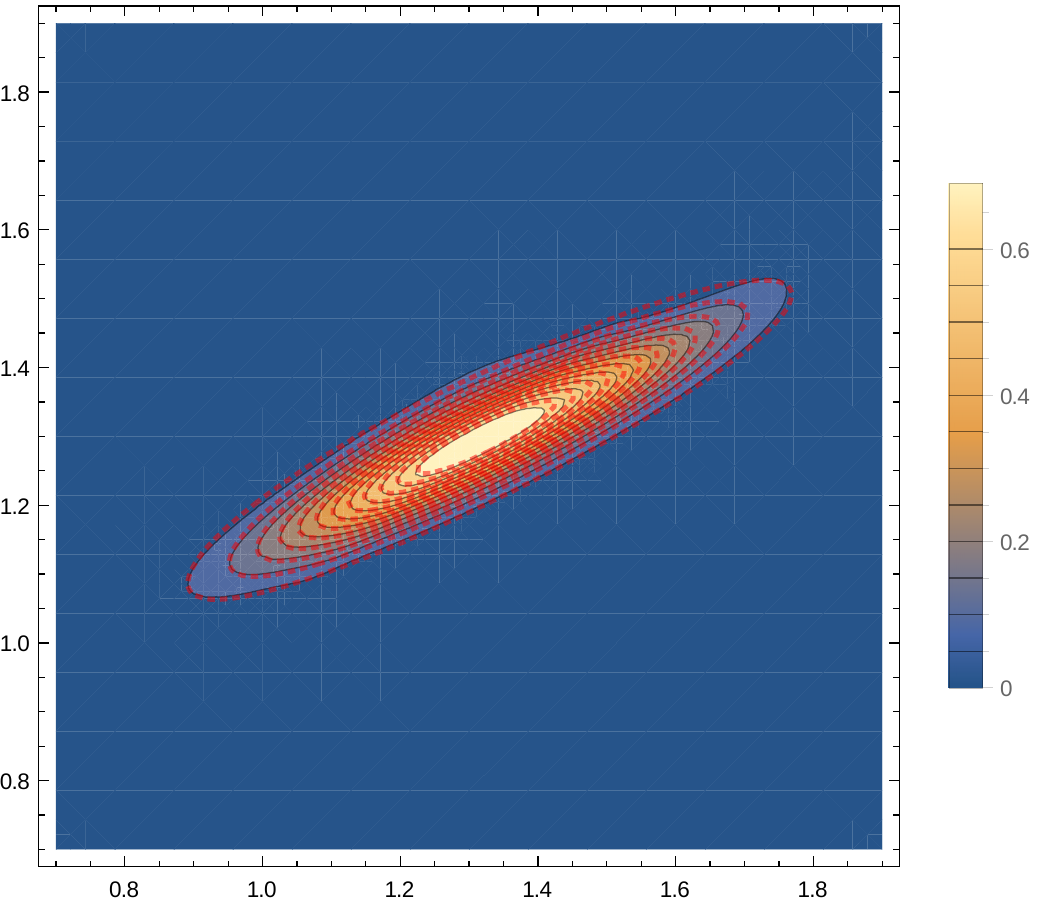}
} \hspace{20pt}
\subfloat[$U_{22}^{\epsilon}$ and the reference solution]{
\includegraphics[height=.35\textwidth]{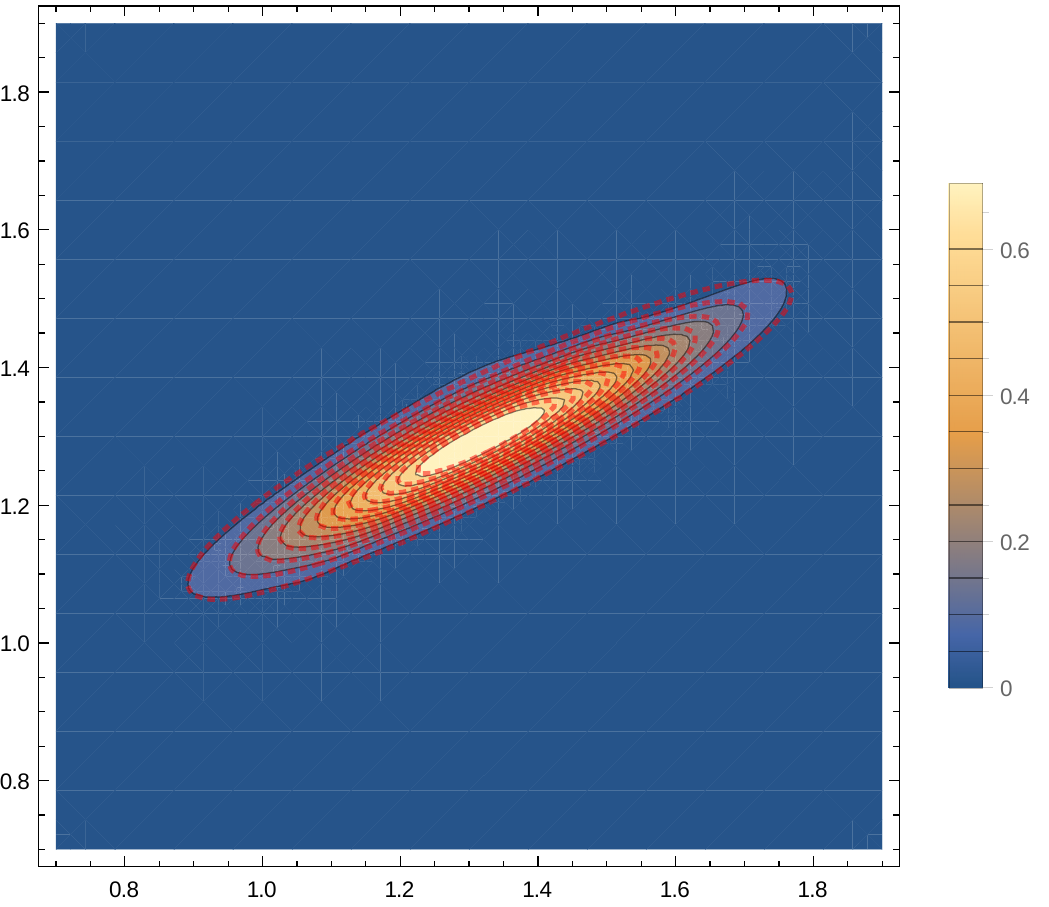}
}
\caption{Comparison between the numerical results (color shading and
black contour lines) and the reference solutions (red dotted contour
lines) for $\epsilon = 1/64$.}
\label{fig:eps_1d64}
\end{figure}

\begin{figure}[!ht]
\centering
\includegraphics[height=.35\textwidth]{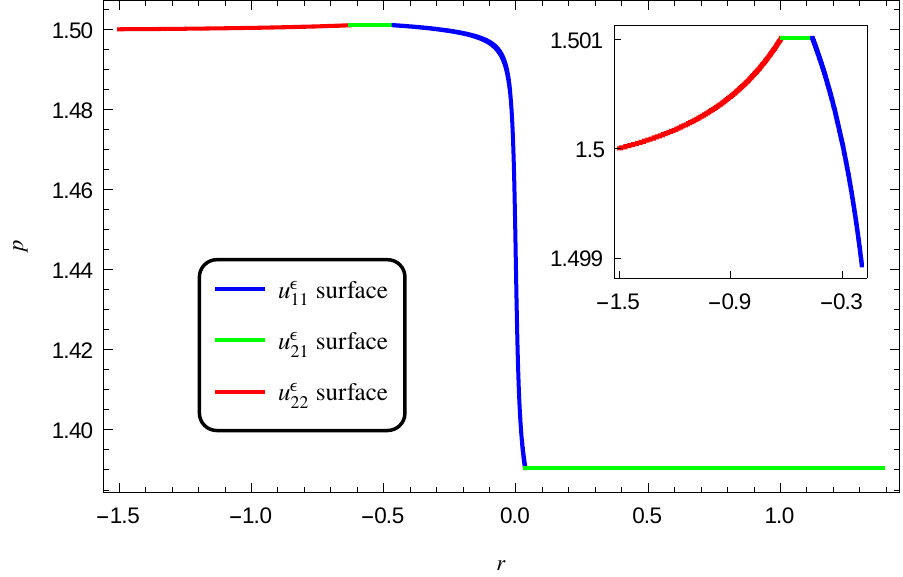}
\caption{A sample trajectory with three surface hops}
\label{fig:traj}
\end{figure}

Due to the space limitations, we only plot figures for
$u_{11}^{\epsilon}$ and $u_{22}^{\epsilon}$. Figure
\ref{fig:eps_1d16}--\ref{fig:eps_1d64} show the comparison with the
reference solutions. In all cases, 10,000 trajectories are used and
the figures show the solutions at $t = 2$. We again observe that
better results are obtained for small $\epsilon$. Figure
\ref{fig:traj} gives a sample of the trajectory in the phase space
$(r,p)$, which is extracted from a test run for $\epsilon = 1/32$.
When $r < 0$, both $E_1(r)$ and $E_2(r)$ are flat, and thus there is
very little change in the momentum $p$. By zooming in, we can still
see that the momentum is increasing slowly until the hop occurs. The
first hop brings the trajectory onto the $u_{21}^{\epsilon}$ surface,
where the momentum does not change since $E_1(r) + E_2(r)$ is zero.
The second hop almost occurs at the point with highest hopping
probability (see Figure \ref{fig:ECR_E}). Since $E_1(r)$ is
increasing, the momentum starts to drop at this hop. When crossing the
point $r = 0$, the momentum has the fastest decrease. The last hop
occurs near $r = 0$, and the momentum again stays as a constant from
then on. Later, since $r$ keeps increasing, the hopping probability
goes lower, no more hops exist in the trajectory.

\subsubsection{Single crossing}
\begin{figure}[!ht]
\centering
\subfloat[Potential energies]{
\includegraphics[height=.3\textwidth]{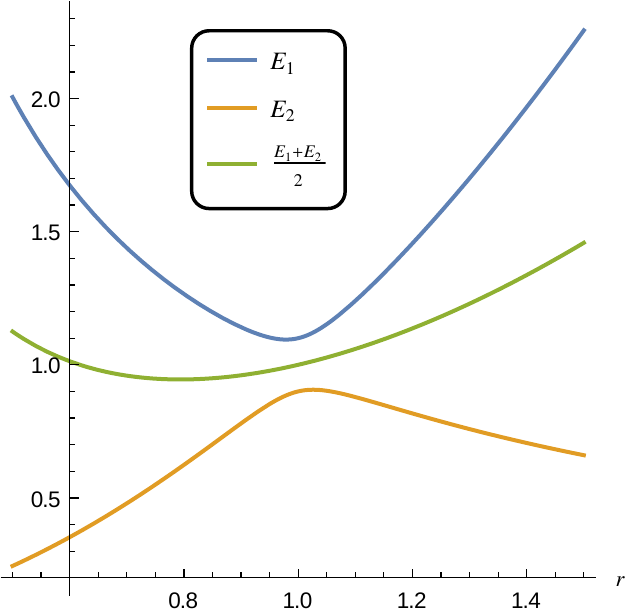}
} \hspace*{50pt}
\subfloat[Nonadiabatic coupling]{
\includegraphics[height=.3\textwidth]{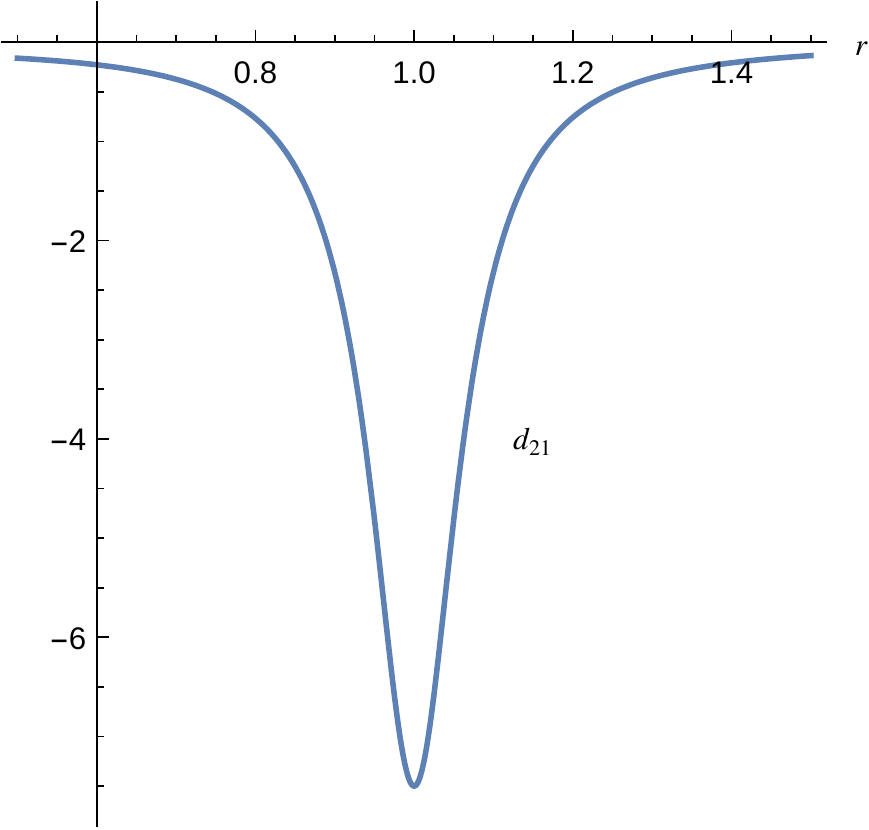}
}
\caption{Potential energies and nonadiabatic coupling for the
single-crossing example}
\label{fig:SC_E}
\end{figure}

This example is a repetition of the example in \cite{Schutte1999,
Horenko2002} with our method. We again consider the one-dimensional
case $x, p \in \mathbb{R}$, and the diabatic potential energy matrix
is
\begin{equation}
V(r) = \begin{pmatrix}
  r^2 & 1/10 \\ 1/10 & 1/r
\end{pmatrix}.
\end{equation}
Thus we have the potential energy $E_{1,2}(r)$, effective potential
energy $\frac{1}{2} [E_1(r) + E_2(r)]$ and the nonadiabatic coupling
coefficient $d_{21}(r)$ as in Figure \ref{fig:SC_E}. The initial value
is
\begin{equation}
u_{11}^{\epsilon}(0,r,p) = \sqrt{\frac{2}{\pi \epsilon}} \exp \left(
  -\frac{(r - r_0)^2}{2\epsilon} - \frac{2(p - p_0)^2}{\epsilon}
\right), \qquad u_{21}^{\epsilon}(0,r,p) = u_{22}^{\epsilon}(0,r,p)=0.
\end{equation}
Following the method in \cite{Yin2013}, we write the initial value as
\begin{equation}
\begin{split}
u_{11}^{\epsilon}(0,r,p) & =
  \int_{\mathbb{R}} \int_{\mathbb{R}} \int_{\mathbb{R}} \int_{\mathbb{R}}
    f(r_1', p_1', r_2', p_2') A(r_1', p_1', r_2', p_2') \times \\
& \qquad \qquad \exp \left(
  -\frac{r^2 + p^2}{\epsilon}
  -\frac{\mathrm{i} [p(r_1'-r_2')-r(p_1'-p_2')]}{\epsilon}
  +\frac{\mathrm{i} (p_1'-p_2')(r_1'+r_2')}{2\epsilon}
\right) \,\mathrm{d}r_1' \,\mathrm{d}p_1'
  \,\mathrm{d}r_2' \,\mathrm{d}p_2',
\end{split}
\end{equation}
where
\begin{align}
\label{eq:rp}
f(r_1', p_1', r_2', p_2') &= \frac{1}{18(\pi\epsilon)^2} \exp \left(
  -\frac{(r_1'-r_0)^2 + (r_2'-r_0)^2}{6\epsilon}
  -\frac{(p_1'-p_0)^2 + (p_2'-p_0)^2}{3\epsilon}
\right), \\
A(r_1', p_1', r_2', p_2') &= \frac{6}{\sqrt{\pi\epsilon}} \exp \left(
  \frac{\mathrm{i}(p_1' - p_0)(2r_1' + r_0)}{3\epsilon}
  -\frac{\mathrm{i}(p_2' - p_0)(2r_2' + r_0)}{3\epsilon}
\right).
\end{align}
Thus the initial condition can be set by drawing $r_1'$, $p_1'$,
$r_2'$ and $p_2'$ from \eqref{eq:rp}.

\begin{figure}[!ht]
\centering
\subfloat[$U_{11}^{\epsilon}$ and the reference solution]{
\includegraphics[height=.33\textwidth]{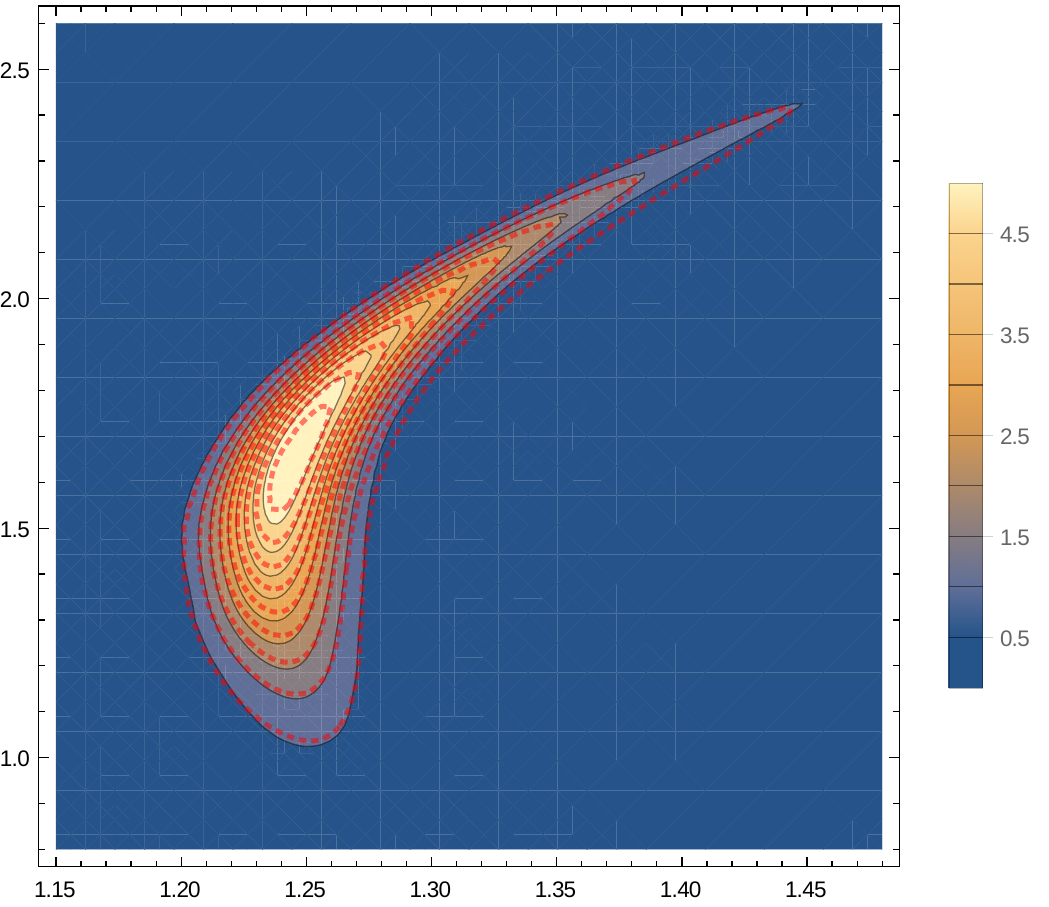}
} \hspace{20pt}
\subfloat[$U_{22}^{\epsilon}$ and the reference solution]{
\includegraphics[height=.33\textwidth]{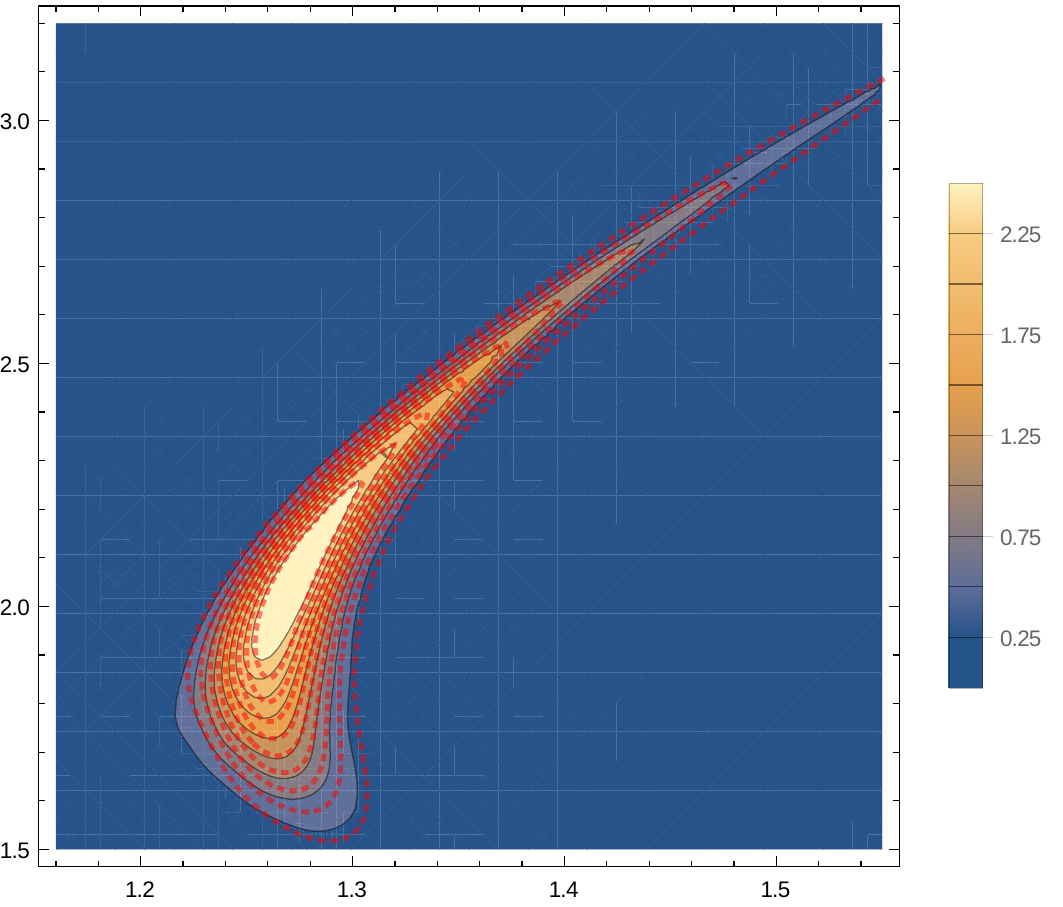}
}
\caption{Comparison between the numerical results (color shading and
black contour lines) and the reference solutions (red dotted contour
lines) for the single-crossing example.}
\label{fig:SC_sol}
\end{figure}

\begin{figure}[!ht]
\centering
\includegraphics[height=.33\textwidth]{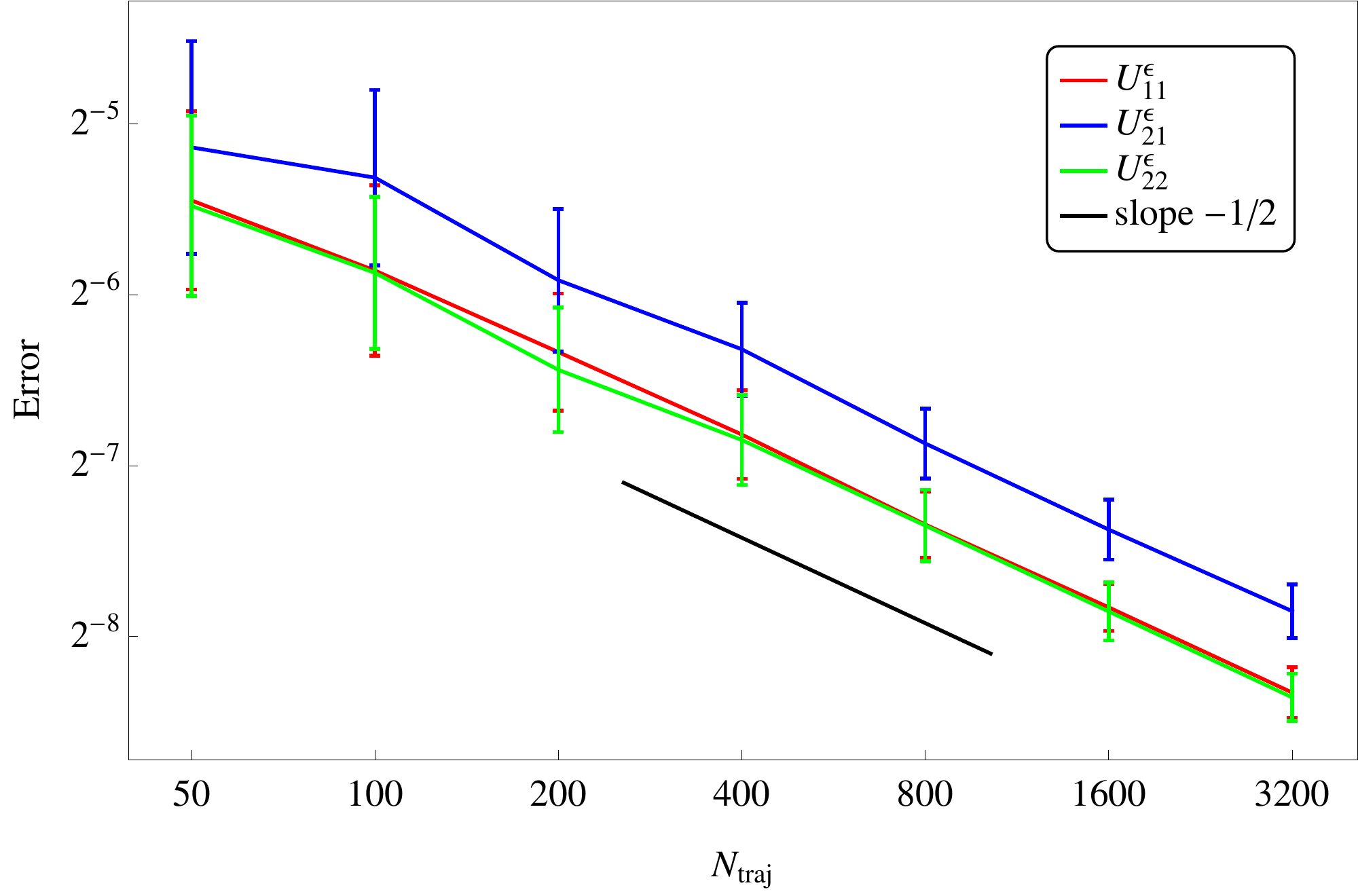}
\caption{$L^2$-error for different numbers of trajectories}
\label{fig:SC_error}
\end{figure}

In our experiment, we follow \cite{Schutte1999, Horenko2002} and
choose $r_0 = 0.4$, $p_0 = 1$ and $\epsilon = 0.01$. Figure
\ref{fig:SC_sol} gives the numerical solution of $U_{11}^{\epsilon}$
and $U_{22}^{\epsilon}$ using 2,000,000 trajectories at $t = 0.5$.
Unlike the last example, the solution no longer has a Gaussian shape
after evolution. Our method can still well capture the structure of
the solution both qualitatively and quantitatively. Using this
as a reference solution, we measure the error of the numerical results
for different numbers of trajectories, which are plotted in Figure
\ref{fig:SC_error}. The errors are again obtained by averaging of 100
runs in each case. The half-order of convergence is still observed.

\section{Summary} \label{sec:summary} We presented a stochastic
numerical method for a type of linear hyperbolic equations with stiff
non-coupling source terms and non-stiff coupling source terms. The
algorithm combines the Gaussian beam method and the surface hopping
method. Similar as \cite{Lu2016}, this work establishes a surface
hopping method based on a solid mathematical foundation, and therefore
gives a better understanding of the surface hopping in a general
setting. As an algorithm, it is very efficient for initial data with
small variance, and works well with the parallelism. As elementary 
application on the quantum-classical dynamics has been shown, and it
can be expected that the algorithm is capable of solving systems with
a large number of particles (corresponding to higher spatial dimension),
which will be explored in the future.

\appendix
\section{Derivation of the expectation} \label{sec:appendix}
Here we present the derivation of the formula \eqref{eq:exp_series}.
Define $F(t;T_1,\cdots,T_K;i_0,i_1,\cdots,i_K)$ as the probability
that $\ell$ jumps $K$ times in $[0,t)$, and the $k$th jump occurs in
$[0, T_k]$ and is from $i_{k-1}$ to $i_k$ (therefore $i_{k-1} \neq
i_k$ for every $k = 1,\cdots,K$). Obviously $F(t;;) =
\mathbb{P}(J_{0,t} = 0)$, and we have
\begin{equation} \label{eq:recur}
F(t; T_1,\cdots,T_K; i_0,\cdots,i_K) = \int_0^{\min(t,T_K)}
  F(t_K; T_1,\cdots,T_{K-1}; i_0,\cdots,i_{K-1}) \Omega_{i_{K-1}i_K}(t_K)
  \mathbb{P}(J_{t_K,t} = 0)
\,\mathrm{d}t_K,
\end{equation}
where $t_K$ means the time of the $K$th jump. The recursion relation
\eqref{eq:recur} leads to
\begin{equation}
\begin{split}
& F(t; T_1,\cdots,T_K; i_0,\cdots,i_K) = \\
& \quad \int_0^{\min(t,T_K)} \int_0^{\min(t_K,T_{k-1})}
  \cdots \int_0^{\min(t_2,T_1)}
  \left( \prod_{k=1}^K \Omega_{i_{k-1} i_k}(t_k) \right)
  \left( \prod_{k=0}^K \mathbb{P}(J_{t_k,t_{k+1}} = 0) \right)
\,\mathrm{d}t_1 \cdots \,\mathrm{d}t_{K-1} \,\mathrm{d}t_K,
\end{split}
\end{equation}
where $t_0 = 0$ and $t_{K+1} = t$. This formula clearly shows that
given $t$ and $i_0,\cdots,i_K$, the conditional probability density
that the $k$th jump occurs at $t_k$ is
\begin{equation}
\rho(t_1,\cdots,t_K \mid t; i_0,\cdots,i_K) = \begin{cases}
  0, \quad \text{if } t_{k+1} < t_k \text{ for any } k=0,\cdots,K, \\
\displaystyle
\frac{1}{F(t; i_0,\cdots,i_K)} \left(
  \prod_{k=1}^K \Omega_{i_{k-1} i_k}(t_k)
\right) \exp \left(
  \int_0^t \Omega_{l_s l_s}(s) \,\mathrm{d}s
\right), \quad \text{otherwise},
\end{cases}
\end{equation}
where $F(t; i_0,\cdots,i_K) := F(t; t,\cdots,t; i_0, \cdots, i_K)$ and
$l_s$ is the jump process defined in \eqref{eq:jump}. Consequently, the
expectation \eqref{eq:exp} can be written as
\begin{equation}
\begin{split}
& \mathbb{E}[\mathcal{F}_i(\ell;t,x)] \\
={} & \sum_{K=0}^{+\infty}
  \sum_{\substack{i_1=1\\i_1\neq i_0}}^n \cdots
  \sum_{\substack{i_K=1\\i_K\neq i_{K-1}}}^n
  \int_0^{+\infty} \!\! \cdots \! \int_0^{+\infty}
    \mathcal{F}_i(\ell;t,x) \rho(t_1,\cdots,t_K \mid t; i_0,\cdots,i_K)
    F(t; i_0, \cdots, i_K)
  \,\mathrm{d}t_1 \cdots \,\mathrm{d}t_K \\
={} & \sum_{K=0}^{+\infty}
  \sum_{\substack{i_1=1\\i_1\neq i_0}}^n \cdots
  \sum_{\substack{i_K=1\\i_K\neq i_{K-1}}}^n
  \int_0^t \!\! \int_0^{t_K} \!\! \cdots \! \int_0^{t_2}
    \mathcal{F}_i(\ell;t,x)
    \left( \prod_{k=1}^K \Omega_{i_k i_{k-1}}(t_k) \right)
    \exp \left( \int_0^t \Omega_{l_s l_s}(s) \,\mathrm{d}s \right)
  \,\mathrm{d}t_1 \cdots \,\mathrm{d}t_{K-1} \,\mathrm{d}t_K.
\end{split}
\end{equation}
Here $\ell$ is the jump process defined by \eqref{eq:jump}. This
result is identical to \eqref{eq:exp_series}.

\bibliographystyle{amsplain}
\bibliography{shgb}
\end{document}